\documentclass[11pt,a4paper,reqno]{amsart}
\usepackage[T1]{fontenc}
\usepackage[utf8]{inputenc}
\usepackage{indentfirst}
\usepackage{lmodern,microtype}
\usepackage[margin=27mm]{geometry}
\usepackage{amssymb,mathtools}
\usepackage{enumitem,needspace}
\setlist[enumerate]{label=\textup{(\arabic*)},leftmargin=2em,itemsep=4pt,topsep=5pt}
\numberwithin{equation}{section}
\usepackage[hidelinks]{hyperref}
\hypersetup{pdftitle={(alpha,beta)-lineability and spaceability of functions with a unique maximum},pdfauthor={Anderson Barbosa, Daniel Pellegrino, Anselmo Raposo Jr.}}
\newtheorem{theorem}{Theorem}[section]
\newtheorem{lemma}[theorem]{Lemma}
\newtheorem{proposition}[theorem]{Proposition}
\newtheorem{corollary}[theorem]{Corollary}
\theoremstyle{remark}
\newtheorem{remark}[theorem]{Remark}
\newtheorem{example}[theorem]{Example}
\theoremstyle{plain}
\newtheorem{maintheorem}{Theorem}

\newcommand{\R}{\mathbb R}
\newcommand{\N}{\mathbb N}
\newcommand{\C}{\mathcal C}
\newcommand{\Chat}{\hat{\mathcal C}}
\newcommand{\cc}{\mathfrak c}
\newcommand{\one}{\mathbf 1}
\DeclareMathOperator{\card}{card}
\DeclareMathOperator{\conv}{conv}
\DeclareMathOperator{\spanop}{span}
\DeclareMathOperator{\Argmax}{Argmax}
\DeclareMathOperator{\Argmin}{Argmin}
\DeclareMathOperator{\interior}{int}
\DeclareMathOperator{\dist}{dist}
\title[Unique maxima and $(\alpha,\beta)$-lineability]{$(\alpha,\beta)$-lineability and spaceability\\ of functions with a unique maximum}
\author{Anderson Barbosa}
\address{Departamento de Ci\^encias Exatas e Tecnologia da Informa\c{c}\~ao,
Universidade Federal Rural do Semi-\'Arido, Rio Grande do Norte, Brazil}
\email{anderson.barbosa@ufersa.edu.br}
\author{Daniel Pellegrino}
\address{Departamento de Matem\'atica, Universidade Federal da Para\'iba,
Jo\~ao Pessoa, Brazil}
\email{daniel.pellegrino@academico.ufpb.br}
\author{Anselmo Raposo Jr.}
\address{Coordena\c{c}\~ao do Curso de Matem\'atica --- Bacharelado,
Universidade Federal do Maranh\~ao, S\~ao Lu\'is, Brazil}
\email{anselmo.junior@ufma.br}
\date{}

\begin{document}
\begin{abstract}
Motivated by the classical work of Gurariy and Quarta on linear spaces
of continuous functions with a unique point of maximum, we determine
the largest possible dimension of such spaces in several settings.
For a connected $n$-manifold $M$ without boundary, with $n\ge2$, the
largest dimension of a subspace of $\C_0(M)$ whose nonzero elements
attain their maximum at exactly one point is $n+1$ if $M$ is
homeomorphic to $S^n$, and $n$ otherwise. We also establish the sharp
bound $n+1$ in $\C(\mathbb R^n)$ and obtain results on extensions of
prescribed finite-dimensional subspaces.
\end{abstract}
\maketitle

\section{Introduction}

Gurariy and Quarta \cite{GQ} initiated in 2004 the study of linear
subspaces of continuous functions whose nonzero elements attain their
maximum at exactly one point. On a compact nondegenerate interval
the dimension is at most one
\cite[Theorem 6]{GQ}, whereas two-dimensional examples exist on
half-open intervals and on $\R$
\cite[Proposition 8 and Theorem 9]{GQ}.
For functions on $\R$ vanishing at infinity, they also established
the sharp bound two \cite[Theorem 16]{GQ}.

In \cite{BCFPS}, constructions involving embedded spheres were used
to obtain finite-dimensional examples on more general domains. The
same work established the upper bound $m$ on compact subsets of
$\R^m$, for every positive integer $m$.
Corollary~\ref{cor:lp-realization} below generalizes the compact metric
example of \cite[Example 4.1]{BCFPS} and identifies a family of
maximally spaceable examples.
In 2020, Bernal-Gonz\'alez, Cabana-M\'endez, Mu\~noz-Fern\'andez and
Seoane-Sep\'ulveda settled the original question on $\R$ by proving
that two is the largest possible dimension in $\C(\R)$ and in
$\C([a,b))$, for $a<b$ \cite[Theorem 3.22]{BCMS}.
Their Section 4 establishes the same sharp bound on $\R$ when
every nonzero function has exactly $m$ maximizing points, for a
fixed integer $m\ge2$.

For functions vanishing at infinity on $\R^n$, with $n\ge2$,
the results in \cite[Theorems 6 and 7]{AB} give admissible subspaces
of dimension $n$ and exclude dimension $n+2$.
Whether dimension $n+1$ can occur remained open
\cite[Conjecture 1]{AB}. The same question is also recorded in
\cite[Conjecture 1.11]{Barbosa}.
Aron and Seoane-Sep\'ulveda's 2025 account
\cite[Example 2.1]{AS} recalls the unique-maximum bounds on compact
intervals and on $\R$.
General background on lineability and spaceability can be found in
the monograph \cite{ABPS} and the survey \cite{BGPS}.

We examine both the possible dimensions and the extension of prescribed
subspaces. All vector spaces and functions are real. Dimensions of vector spaces
are Hamel dimensions; covering dimension of a topological space is
written $\dim_{\mathrm{top}}$ and defined in Section~\ref{sec:proofC}. Write
$\N=\{1,2,\ldots\}$, $\aleph_0=\card(\N)$ and $\cc=\card(\R)$.
For $d\in\N$, $\ell_2^d$ denotes $\R^d$ with its Euclidean norm.
For a real function $f:X\to\R$ on a nonempty set, put
\[
 \Argmax_X f=\{x\in X:f(x)=\sup_{z\in X}f(z)\},\qquad
 \Argmin_X f=\Argmax_X(-f).
\]
The supremum is taken in the extended real line. Thus $\Argmax_X f$
is empty when $f$ is unbounded above or its supremum is not attained;
the condition $|\Argmax_X f|=1$ means that $f$ attains its maximum
at exactly one point.
Following the notation of \cite{AB}, let $\C(X)$ be the space of
continuous real functions on $X$, with
the compact-open topology given by the seminorms
$f\mapsto\sup_{x\in H}|f(x)|$ for nonempty compact $H\subset X$.
Let $\C_b(X)\subset \C(X)$ consist of the bounded functions, with
norm $\|f\|_\infty=\sup_{x\in X}|f(x)|$, and define
\[
 \Chat(X)=\{f\in \C_b(X):|\Argmax_Xf|=1\}.
\]
For a Hausdorff space $X$, let $\C_0(X)$
be the space of continuous functions $f:X\to\R$ such that
$\{x:|f(x)|\ge\varepsilon\}$ is compact for every $\varepsilon>0$,
with the supremum norm. This definition also applies when $X$ is not
locally compact; in the locally compact Hausdorff case it is the usual
Banach space of functions vanishing at infinity. Define
\[
 \Chat_0(X)=\{f\in \C_0(X):|\Argmax_Xf|=1\}.
\]
For compact Hausdorff $K$, we have
$\C_0(K)=\C(K)$ and $\Chat_0(K)=\Chat(K)$.
A linear subspace of $\C(X)$ whose nonzero members attain unique
maxima is called \emph{admissible}. Although a continuous function with
a unique maximum need not be bounded below, every admissible subspace
consists of bounded functions, since it contains both $f$ and $-f$,
which attain finite maxima. Its nonzero members also attain unique
minima. Thus every admissible subspace of $\C(X)$
lies in $\Chat(X)\cup\{0\}$. We consider these subspaces in
$\C(X)$ and in $\C_b(X)$, with the respective topologies specified above.

The notion of $(\alpha,\beta)$-lineability was introduced in
\cite[Definition 1.1]{FPT}. A subset $A$ of a real vector space $V$
is $\alpha$-lineable if $A\cup\{0\}$ contains an
$\alpha$-dimensional subspace. For cardinals $\alpha\le\beta$,
it is $(\alpha,\beta)$-lineable if it is $\alpha$-lineable and
every subspace $E\subset A\cup\{0\}$ of dimension $\alpha$ is
contained in a subspace $W\subset A\cup\{0\}$ of dimension $\beta$.
In a topological vector space, $(\alpha,\beta)$-spaceability
requires that $A$ be $\alpha$-lineable and a \emph{closed} extension
$W$ with these properties. The prescribed $E$ need not be closed.
We include pairs with $\alpha=\beta$, as in \cite{FPRR}. A set is $\alpha$-spaceable
if it contains, together with zero, a closed $\alpha$-dimensional subspace;
it is maximally spaceable when $\alpha$ is the dimension of the ambient space.
For normed spaces $E_1,E_2$ and a bounded linear map
$T:E_1\to E_2$, write
$\|T\|=\sup_{\|f\|_{E_1}\le1}\|Tf\|_{E_2}$.
In 2024, results on $(\alpha,\beta)$-lineability and
$(\alpha,\beta)$-spaceability for complements of vector
subspaces were obtained in \cite{ABcriterion}, while general
obstructions to closed extensions were established in
\cite[Theorem 2.1]{FPRR}.
A related obstruction appears in \cite[Theorem 9]{AB}; in Banach
spaces its hypotheses require infinite prescribed dimension. The finite
prescribed dimensions considered below require a different argument.
In particular, failure of $(\alpha,\beta)$-lineability does not mean
that no $\alpha$-dimensional subspace extends: one nonextendible
$\alpha$-dimensional subspace is enough to disprove the property.

Recall that $X$ is countably compact if every countable open cover has
a finite subcover. Every continuous map from such a space into a
metrizable space has compact image. In particular, $\C(X)=\C_b(X)$
for countably compact $X$; this does not imply $\C_0(X)=\C(X)$.
For a noncompact $X$, the constant function one does not vanish at infinity.

A finite-dimensional admissible subspace can be perturbed, by an
arbitrarily small isometric change preserving every maximizing and
minimizing point, into an admissible subspace that admits no proper
admissible extension. The precise statement and proof appear in
Theorem~\ref{lem:perturbation}. Here \emph{maximal by inclusion} means
that no proper admissible extension exists, not that the subspace has
the largest possible dimension.

The perturbation theorem gives the following
classification whenever the prescribed dimension is positive and finite.

\Needspace{8\baselineskip}
\begin{maintheorem}\label{thm:A}
Let $X$ be locally compact Hausdorff with at least two points.
For cardinals $1\le\alpha<\aleph_0$ and $\beta\ge\alpha$,
\[
 \Chat_0(X)\text{ is }(\alpha,\beta)\text{-lineable}
 \quad\Longleftrightarrow\quad
 \alpha=\beta\ \text{ and }\ \Chat_0(X)\text{ is }\alpha\text{-lineable}.
\]
The same equivalence holds for $(\alpha,\beta)$-spaceability.
If instead $X$ is countably compact Hausdorff, the same equivalences
hold with $\Chat(X)$ in place of $\Chat_0(X)$, for spaceability in
either the supremum norm or the compact-open topology.
\end{maintheorem}

On a connected $n$-manifold $M$ with $n\ge2$, the largest admissible
dimension in $\C_0(M)$ has a topological description.
By an $n$-manifold we mean a nonempty Hausdorff second-countable space locally
homeomorphic to $\R^n$, without boundary. Write $M\cong N$ for
homeomorphism and
$S^n=\{a\in\R^{n+1}:\sum_{j=1}^{n+1}a_j^2=1\}$.
\Needspace{10\baselineskip}
\begin{maintheorem}\label{thm:B}
Let $M$ be a connected $n$-manifold, with $n\ge2$, and let
$1\le\alpha\le\beta$ be cardinals.
\begin{enumerate}
\item If $M\cong S^n$, then $\Chat_0(M)$ is
$(\alpha,\beta)$-lineable if and only if $\alpha=\beta\le n+1$.
\item If $M\not\cong S^n$, then $\Chat_0(M)$ is
$(\alpha,\beta)$-lineable if and only if $\alpha=\beta\le n$.
\end{enumerate}
The same equivalences hold for $(\alpha,\beta)$-spaceability.
\end{maintheorem}

In particular, $\Chat_0(\R^n)$ is $n$-lineable but not
$(n+1)$-lineable for every $n\ge2$, resolving
\cite[Conjecture 1]{AB}.
The case $n=1$ is an essential exception: $\Chat_0(\R)$ has largest
admissible dimension two \cite[Theorem 16]{GQ}. Together with
Theorem~\ref{thm:A}, this gives $\alpha=\beta\le2$ on $\R$.

For a topological space $X$, a \emph{zero-set} is
$u^{-1}(\{0\})$ for some continuous $u:X\to[0,1]$.
Continuous real functions separate points if, whenever $x\ne y$,
there is a continuous $f:X\to\mathbb R$ with $f(x)\ne f(y)$.
Proposition~\ref{prop:countably-compact-sphere} characterizes finite
admissible dimensions on countably compact spaces with this separation
property by zero-sets homeomorphic to spheres. Combined with
Theorem~\ref{thm:A} and the infinite-dimensional obstruction in
Lemma~\ref{lem:constantshift}, it gives the spaceability criterion
in Corollary~\ref{cor:countable-classification}.
The sphere must be a zero-set; for compact metrizable domains,
any embedded sphere has this property.

Corollary~\ref{cor:transfer} transfers the finite-dimensional
nonextension result from $\C_0(X)$ to $\C(X)$. On Euclidean domains,
we also determine the largest admissible dimension without imposing
vanishing at infinity.
The upper bound extends to every $\sigma$-compact metrizable
space of finite covering dimension; on Euclidean domains it is attained.

\Needspace{10\baselineskip}
\begin{maintheorem}\label{thm:C}
For every integer $n\ge1$, the set $\Chat(\mathbb R^n)$ is
$(n+1)$-lineable but not $(n+2)$-lineable. In fact, it contains,
together with zero, a linear isometric copy of $\ell_2^{n+1}$.
For all cardinals $1\le\alpha\le\beta$, each of the following
properties holds if and only if $1\le\alpha=\beta\le n+1$:
\begin{enumerate}
\item $(\alpha,\beta)$-lineability of $\Chat(\mathbb R^n)$;
\item $(\alpha,\beta)$-spaceability in $\C_b(\mathbb R^n)$ with the
supremum norm;
\item $(\alpha,\beta)$-spaceability in $\C(\mathbb R^n)$ with the
compact-open topology.
\end{enumerate}
\end{maintheorem}

For $n\ge2$, the largest admissible dimensions in
$\C_0(\R^n)$ and $\C(\R^n)$ are $n$ and $n+1$, respectively.
For $n=1$, both dimensions equal two.

\section{Proof of Theorem A}\label{sec:proofA}

A proper admissible extension forces a continuous section of the
evaluation map. A small radial deformation preserves all extrema
while preventing such a section.

\subsection{Convex sets and evaluation maps}\label{sec:universalauxiliary}

The convex-hull arguments are related to
\cite[Propositions 3.1 and 3.3, Lemma 3.5 and Theorem 3.8]{BCMS} and
\cite[Lemmas 3.2--3.5]{BCFPS}.
Continuity of the maximizing-point map on an interval already appears
in \cite[Lemma 5]{GQ}.
We use the Euclidean inner product
$\langle a,u\rangle=\sum_{j=1}^d a_ju_j$ and norm
$\|a\|=(\sum_{j=1}^d a_j^2)^{1/2}$ on $\mathbb R^d$.
For $C\subset\mathbb R^d$, its Euclidean interior and boundary are
written $\interior C$ and $\partial C$. We refer to \cite{Rockafellar}
for the standard convex-geometric facts used below, including
Carath\'eodory's theorem, supporting hyperplanes and strict separation.

\begin{lemma}\label{lem:evaluation}
Let $Y$ be a nonempty topological space, let $d\in\mathbb N$,
and let $\Psi:Y\to\mathbb R^d$ be continuous with compact image.
This image condition holds, in particular, when $Y$ is countably compact.
Let $C=\conv\Psi(Y)$,
where $\conv$ denotes the set of finite convex combinations.
For $a\in\mathbb R^d\setminus\{0\}$, define
\[
 q_a:Y\to\mathbb R,\qquad q_a(x)=\langle a,\Psi(x)\rangle,
\]
and put
\[
 H_C(a)=\max_{u\in C}\langle a,u\rangle,\qquad
 C_a=\{u\in C:\langle a,u\rangle=H_C(a)\}.
\]
Then the following statements hold.
\begin{enumerate}
\item The set $C$ is compact, and
\begin{equation}\label{eq:evaluationface}
 H_C(a)=\max_{x\in Y}q_a(x),\qquad
 C_a=\conv\Psi(\Argmax_Yq_a).
\end{equation}
\item If $\Argmax_Yq_a=\{x_a\}$, then
$C_a=\{\Psi(x_a)\}$ and
$\Psi^{-1}(\{\Psi(x_a)\})=\{x_a\}$.
\end{enumerate}
\end{lemma}

\begin{proof}
	\textup{(1)} By Carath\'eodory's theorem, every point of $C$ is a
	convex combination of at most $d+1$ points of $\Psi(Y)$.
	Thus $C$ is the image of the compact space
	\[
	\left\{(\lambda_1,\ldots,\lambda_{d+1})\in[0,1]^{d+1}:
	\sum_{j=1}^{d+1}\lambda_j=1\right\}\times \Psi(Y)^{d+1}
	\]
	under the continuous map
	$((\lambda_j)_{j=1}^{d+1},(x_j)_{j=1}^{d+1})
	\mapsto\sum_{j=1}^{d+1}\lambda_j\Psi(x_j)$.
	In particular, $C$ is compact and the maximum $H_C(a)$ is attained.
	For any convex combination $u=\sum_{j=1}^r\lambda_j\Psi(x_j)$,
	with $r\in\mathbb N$, $\lambda_j\ge0$ and
	$\sum_{j=1}^r\lambda_j=1$, we have
	\[
	\langle a,u\rangle
	=\sum_{j=1}^r\lambda_jq_a(x_j)
	\le\sum_{j=1}^r\lambda_j\max_Yq_a
	=\max_Yq_a.
	\]
	Since $\Psi(Y)\subset C$, this proves the first identity in
	\eqref{eq:evaluationface}.
	
	Write $H=H_C(a)$. If $u\in C_a$, the same representation gives
	\[
	0=H-\langle a,u\rangle
	=\sum_{j=1}^r\lambda_j\bigl(H-q_a(x_j)\bigr).
	\]
	Every summand is nonnegative. Therefore $q_a(x_j)=H$ whenever
	$\lambda_j>0$, so $u\in\conv\Psi(\Argmax_Yq_a)$.
	Conversely, if every $x_j$ with $\lambda_j>0$ belongs to
	$\Argmax_Yq_a$, then
	$\langle a,u\rangle=\sum_{j=1}^r\lambda_jH=H$.
	This proves $C_a=\conv\Psi(\Argmax_Yq_a)$.
	
	\smallskip
	\noindent\textup{(2)} The equality $C_a=\conv\Psi(\Argmax_Yq_a)$ gives $C_a=\{\Psi(x_a)\}$.
	If $\Psi(x)=\Psi(x_a)$, then
	$q_a(x)=\langle a,\Psi(x_a)\rangle=H$.
	Uniqueness of the maximizer yields $x=x_a$.
\end{proof}

The Minkowski functional provides radial coordinates on a convex body
and a parametrization of its boundary by a sphere.

\begin{lemma}\label{lem:gauge}
Let $C\subset\mathbb R^k$ be compact and convex with nonempty interior,
where $k\in\mathbb N$, and let $c\in\interior C$. Define the
Minkowski functional of $C-c$ by
\[
 p:\mathbb R^k\to[0,\infty),\qquad
 p(y)=\inf\{t>0:y\in t(C-c)\}.
\]
\begin{enumerate}
\item The function $p$ is continuous, $p(y)=0$ exactly when $y=0$,
and $p(ty)=tp(y)$ for $t\ge0$. Moreover,
\begin{equation}\label{eq:gaugesets}
 \begin{aligned}
 C-c&=\{y:p(y)\le1\},\\
 \interior(C-c)&=\{y:p(y)<1\},\qquad
 \partial(C-c)=\{y:p(y)=1\}.
 \end{aligned}
\end{equation}
\item With $S^{k-1}=\{a\in\mathbb R^k:\|a\|=1\}$, the map
\[
 R:S^{k-1}\to\partial C,\qquad R(a)=c+\frac{a}{p(a)},
\]
is a homeomorphism. Its inverse is
$R^{-1}(u)=(u-c)/\|u-c\|$.
\end{enumerate}
\end{lemma}

\begin{proof}
\textup{(1)} Put $Q=C-c$. Since $0\in\interior Q$ and $Q$ is compact,
there are $0<\rho\le R_0$ such that
\[
 \{y:\|y\|\le\rho\}\subset Q
 \subset\{y:\|y\|\le R_0\}.
\]
If $y\in tQ$, then $\|y\|\le tR_0$. Conversely, for $y\ne0$,
the vector $\rho y/\|y\|$ belongs to $Q$. Consequently
\begin{equation}\label{eq:gaugebounds}
 \frac{\|y\|}{R_0}\le p(y)\le\frac{\|y\|}{\rho}
 \qquad(y\in\mathbb R^k).
\end{equation}
For $t>0$, substituting $s=tr$ in the defining infimum gives
\[
 p(ty)=\inf\{tr:r>0,
                     ty\in trQ\}
       =t\inf\{r>0:y\in rQ\}=tp(y).
\]
For $t=0$, the same identity follows from $p(0)=0$.

Convexity and $0\in Q$ imply $sQ\subset tQ$ when $0\le s\le t$.
Hence $s>p(y)$ implies $y\in sQ$: there exists $r\in(0,s)$ such that $y\in rQ$,
and $rQ\subset sQ$.
If also $t>p(z)$, convexity gives
\[
 \frac{y+z}{s+t}
 =\frac{s}{s+t}\frac{y}{s}
  +\frac{t}{s+t}\frac{z}{t}\in Q.
\]
Thus $p(y+z)\le s+t$. Letting $s\downarrow p(y)$ and
$t\downarrow p(z)$ proves subadditivity. In particular,
\[
 p(y)-p(z)\le p(y-z)\le\frac{\|y-z\|}{\rho},\qquad
 p(z)-p(y)\le p(z-y)\le\frac{\|y-z\|}{\rho}.
\]
Therefore $|p(y)-p(z)|\le\|y-z\|/\rho$, proving continuity.
The lower bound in \eqref{eq:gaugebounds} proves the assertion about zeros.

If $y\in Q$, then $p(y)\le1$. If $p(y)\le1$, then
$y/(1+1/m)\in Q$ for every $m\in\mathbb N$; closedness of $Q$ gives
$y\in Q$ after passage to the limit.
Next, suppose $p(y)<1$. Choose $s\in(p(y),1)$ and write $y=sz$ with
$z\in Q$. For $\|b\|<(1-s)\rho$,
\[
 y+b=sz+(1-s)\frac{b}{1-s}\in Q,
\]
because $z\in Q$ and $b/(1-s)\in Q$.
Thus $y$ lies in the interior of $Q$.
Conversely, if $0\ne y\in\interior Q$, then
$(1+\delta)y\in Q$ for some $\delta>0$, whence
$p(y)\le(1+\delta)^{-1}<1$ by homogeneity.
For $y=0$, this follows from $p(0)=0$.
Since $Q$ is closed, $\partial Q=Q\setminus\interior Q$;
this proves all of \eqref{eq:gaugesets}.

\smallskip
\noindent\textup{(2)} For $a\in S^{k-1}$, the lower bound in
\eqref{eq:gaugebounds} gives $p(a)>0$, and homogeneity gives
$p(a/p(a))=1$. Hence $R$ is well defined and continuous.
For $u\in\partial C$, let $a=(u-c)/\|u-c\|$.
Since $p(u-c)=1$, we obtain
\[
 p(a)=\frac{1}{\|u-c\|},\qquad
 R(a)=c+\frac{(u-c)/\|u-c\|}{1/\|u-c\|}=u.
\]
Also $(R(a)-c)/\|R(a)-c\|=a$, since $p(a)>0$.
The displayed inverse is continuous, so $R$ is a homeomorphism.
\end{proof}

For a finite-dimensional admissible subspace of $\C_0(X)$,
extension by zero adds an evaluation vector at infinity. The argument
below distinguishes whether this vector is an interior or boundary
point of the convex hull; the boundary case may produce a segment.

\begin{lemma}\label{lem:basis}
Let $X$ be locally compact Hausdorff with at least two points, and let
$f_1,\ldots,f_k$ be a basis of a subspace
$E\subset\Chat_0(X)\cup\{0\}$, where $k\in\mathbb N$.
Take $L=X$ when $X$ is compact, and otherwise take
$L=X\cup\{\infty\}$, the one-point compactification,
whose neighborhoods of $\infty$ are the sets
$\{\infty\}\cup(X\setminus H)$ with $H\subset X$ compact.
Extend each $f_j$ by zero at $\infty$. Define
\[
 F:L\to\mathbb R^k,\qquad
 F(x)=(f_1(x),\ldots,f_k(x)),\qquad C=\conv F(L).
\]
For $a\ne0$, write $q_a(x)=\langle a,F(x)\rangle$ on $L$, and use
$H_C(a)$ and $C_a$ as in Lemma~\ref{lem:evaluation}.
\begin{enumerate}
\item The set $C$ is compact with nonempty interior, and
$H_C(a)=\max_Xq_a$. If $X$ is noncompact, then $0\in C$ and
$H_C(a)\ge0$ for all $a\ne0$.
\item If $X$ is compact, or if $X$ is noncompact and
$0\in\interior C$, then, for $Z=F^{-1}(\partial C)$, the restriction
$F|_Z:Z\to\partial C$ is a homeomorphism and $Z\subset X$.
\item Suppose $X$ is noncompact and $0\in\partial C$.
There is a point $x_0\in X$ which is the unique zero of every
nonzero nonpositive function in $E$. Put
$A=[0,F(x_0)]=\{tF(x_0):0\le t\le1\}$.
Then $C_a=A$ whenever $H_C(a)=0$,
$F^{-1}(A)=\{x_0,\infty\}$, and every point of
$\partial C\setminus A$ has exactly one preimage, lying in $X$.
\end{enumerate}
\end{lemma}

\begin{proof}
\textup{(1)}
The space $L$ is compact Hausdorff, and the coordinate functions
of $F$ are continuous on $L$. In the noncompact case, their
continuity at $\infty$ follows from $f_j\in\C_0(X)$ and the
definition $f_j(\infty)=0$. Thus $F$ is continuous, and
Lemma~\ref{lem:evaluation}\textup{(1)} shows that
$C=\conv F(L)$ is compact.

We prove that $\interior C\ne\varnothing$.
If the affine span of $C$ were all of $\mathbb R^k$, then $C$
would contain $k+1$ affinely independent points. By convexity,
it would contain their convex hull, which has nonempty
interior in $\mathbb R^k$. Consequently, if $C$ had empty
interior, its affine span would be a proper affine subspace.
There would then exist $a\in\mathbb R^k\setminus\{0\}$ and
$b\in\mathbb R$ such that
\begin{equation}\label{eq:affinerelation}
\sum_{j=1}^k a_jf_j(x)=b
\qquad(x\in L).
\end{equation}
If $b=0$, restriction to $X$ contradicts the linear independence
of $f_1,\ldots,f_k$. If $b\ne0$, that restriction is a nonzero
constant belonging to $E$. Since $X$ has at least two points,
such a function does not have a unique maximizer, again a
contradiction. Therefore $\interior C\ne\varnothing$.

Fix $a\in\mathbb R^k\setminus\{0\}$. The restriction $q_a|_X$
is nonzero by linear independence, so it attains its maximum
at a unique point, which we denote by $x_a$.
Lemma~\ref{lem:evaluation}\textup{(1)} gives
$H_C(a)=\max_Lq_a$. If $X$ is compact, then $L=X$, and hence
$H_C(a)=\max_Xq_a$.

Suppose now that $X$ is noncompact. For every $\varepsilon>0$,
the set
\[
K_\varepsilon
=\{x\in X:|q_a(x)|\ge\varepsilon\}
\]
is compact because $q_a|_X\in\C_0(X)$.
It cannot equal the noncompact space $X$. Choosing
$x\in X\setminus K_\varepsilon$ gives
$\sup_Xq_a\ge q_a(x)>-\varepsilon$.
Since $\varepsilon>0$ is arbitrary, $\sup_Xq_a\ge0$.
This supremum is attained at $x_a$. Adding the value
$q_a(\infty)=0$ therefore does not change it, and we obtain
\[
H_C(a)=\max_Lq_a=\max_Xq_a\ge0.
\]
Also, $F(\infty)=0$, so $0\in C$.

In the noncompact case, the same argument gives the following relation
between the maximizing sets on $X$ and on $L$:
\begin{equation}\label{eq:basis-maximizers}
\Argmax_Lq_a=
\begin{cases}
\{x_a\},&H_C(a)>0,\\
\{x_a,\infty\},&H_C(a)=0.
\end{cases}
\end{equation}
Indeed, $x_a$ is the unique maximizer in $X$, and the only
additional point of $L$ is $\infty$, where $q_a$ has value zero.

\smallskip
\noindent
\textup{(2)}
Suppose that $X$ is compact, or that $X$ is noncompact and
$0\in\interior C$.

In the compact case, every $q_a$, with $a\ne0$, has a unique
maximizer on $L=X$. In the noncompact case, choose $\rho>0$
such that $\{u:\|u\|\le\rho\}\subset C$. Then
\[
H_C(a)\ge
\left\langle a,\rho\frac{a}{\|a\|}\right\rangle
=\rho\|a\|>0=q_a(\infty).
\]
Thus \eqref{eq:basis-maximizers} shows that $q_a$ also has a
unique maximizer on $L$, and that this maximizer belongs to $X$.
In either case, Lemma~\ref{lem:evaluation}\textup{(2)} yields
\begin{equation}\label{eq:basis-singleton}
C_a=\{F(x_a)\},
\qquad
F^{-1}(\{F(x_a)\})=\{x_a\}.
\end{equation}

Let $u\in\partial C$. The supporting hyperplane theorem gives
$a\ne0$ such that
$\langle a,v\rangle\le\langle a,u\rangle$ for every $v\in C$.
Equivalently, $u\in C_a$. By \eqref{eq:basis-singleton},
$u=F(x_a)$ and $x_a$ is its unique preimage in $L$.
Consequently, for $Z=F^{-1}(\partial C)$, the restriction
$F|_Z:Z\to\partial C$ is a continuous bijection.

The set $Z$ is closed in compact $L$, because $F$ is continuous
and $\partial C$ is closed. Thus $Z$ is compact.
Since $\partial C$ is Hausdorff, the continuous bijection
$F|_Z$ has a continuous inverse.
Finally, $Z\subset X$: this is immediate when $L=X$, while
in the noncompact case $F(\infty)=0\in\interior C$ implies
$\infty\notin Z$.

\smallskip
\noindent
\textup{(3)}
Suppose that $X$ is noncompact and $0\in\partial C$.
Every nonzero nonpositive function in $E$
has the same unique zero.

A supporting hyperplane at $0$ gives
$a_0\in\mathbb R^k\setminus\{0\}$ such that
$\langle a_0,u\rangle\le0$ for every $u\in C$.
Since $0\in C$, we have $H_C(a_0)=0$. Therefore the set
\[
\mathcal P
=\{f\in E\setminus\{0\}:f(x)\le0\text{ for every }x\in X\}
\]
is nonempty, as it contains $q_{a_0}|_X$.
By part~\textup{(1)}, every $f\in\mathcal P$ has maximum zero.
Since $f$ has a unique maximizer, it has exactly one zero.

Let $f,g\in\mathcal P$. Their sum is nonzero: if $f+g=0$,
then $f=-g$, and the nonpositivity of both functions would
force $f=g=0$. Hence $f+g\in\mathcal P$, so there exists
$x\in X$ such that $(f+g)(x)=0$. As $f(x)\le0$ and $g(x)\le0$,
this equality implies $f(x)=g(x)=0$.
The unique zeros of $f$ and $g$ therefore coincide.
Fixing one member of $\mathcal P$ and comparing it with every
other member gives a point $x_0\in X$ which is the unique zero
of every function in $\mathcal P$.

With $A=[0,F(x_0)]$, we have
\[
C_a=A\quad\text{whenever }H_C(a)=0,
\qquad
F^{-1}(A)=\{x_0,\infty\}.
\]
Suppose $a\ne0$ and $H_C(a)=0$.
Then $q_a|_X\in\mathcal P$, so its unique maximizer is $x_0$.
By \eqref{eq:basis-maximizers},
$\Argmax_Lq_a=\{x_0,\infty\}$.
Applying the second equality in \eqref{eq:evaluationface}
and using $F(\infty)=0$, we obtain
\[
C_a
=\conv F(\{x_0,\infty\})
=\conv\{F(x_0),0\}
=[0,F(x_0)]=A.
\]
In particular, $A\subset\partial C$. Indeed, every point of
$A=C_a$ maximizes the nonzero functional
$u\mapsto\langle a,u\rangle$ on $C$, whereas no interior point
can do so: a sufficiently small displacement in the direction
$a$ would remain in $C$ and increase the functional.

To determine $F^{-1}(A)$, use the fixed direction $a_0$ chosen
above. Since $q_{a_0}(x_0)=0$, the functional
$u\mapsto\langle a_0,u\rangle$ vanishes at both endpoints of
$A$, and therefore throughout $A$.
If $x\in X$ and $F(x)\in A$, then
$q_{a_0}(x)=\langle a_0,F(x)\rangle=0$.
The uniqueness of the zero of $q_{a_0}|_X$ gives $x=x_0$.
Conversely, $F(x_0)\in A$ and $F(\infty)=0\in A$.
Thus
\[
F^{-1}(A)=\{x_0,\infty\}.
\]
In particular,
\[
F(L)\cap A=\{0,F(x_0)\}.
\]
Thus, if $F(x_0)\ne0$, the points $tF(x_0)$ with $0<t<1$
belong to $C$ but not to $F(L)$; they arise as convex
combinations of the vectors $0=F(\infty)$ and $F(x_0)$.

Finally, let $u\in\partial C\setminus A$.
By the supporting hyperplane theorem, choose $a\ne0$ such that
$\langle a,v\rangle\le\langle a,u\rangle$ for every $v\in C$,
so that $u\in C_a$. Part~\textup{(1)} gives $H_C(a)\ge0$.
If $H_C(a)=0$, the equality just proved would give
$C_a=A$, contradicting $u\notin A$.
Hence $H_C(a)>0$.
By \eqref{eq:basis-maximizers}, $q_a$ has exactly one maximizer
on $L$, and this point belongs to $X$.
Lemma~\ref{lem:evaluation}\textup{(2)} now shows that
$u$ has exactly one preimage in $L$, lying in $X$.
\end{proof}

Adding a basis vector adds a coordinate to the evaluation map.
The resulting convex body projects onto the original one; continuous
sections will be obtained from the endpoints of its vertical fibers.

\begin{lemma}\label{lem:vertical}
Let $k\in\mathbb N$ and let $D\subset\mathbb R^k\times\mathbb R$ be
compact and convex with nonempty interior. Let
$\pi:\mathbb R^k\times\mathbb R\to\mathbb R^k$ be the projection
$\pi(u,z)=u$, and put $C=\pi(D)$.
For $u\in\interior C$, define
\[
 z_-(u)=\min\{z:(u,z)\in D\},\qquad
 z_+(u)=\max\{z:(u,z)\in D\}.
\]
\begin{enumerate}
\item The functions $z_-,z_+:\interior C\to\mathbb R$ are continuous,
and $z_-(u)<z_+(u)$ with $(u,z_\pm(u))\in\partial D$.
\item Let $\ell:\mathbb R^k\times\mathbb R\to\mathbb R$ be a nonzero
linear functional, $\ell(u,z)=\langle a,u\rangle+bz$, with
$\ell\le0$ on $D$. Then $\ell(u,z_-(u))<0$ when $b>0$,
$\ell(u,z_+(u))<0$ when $b<0$, and both inequalities hold when $b=0$.
\end{enumerate}
\end{lemma}

\begin{proof}
\textup{(1)} For each $u\in C$, the set $I_u=\{z\in\mathbb R:(u,z)\in D\}$
is nonempty, compact and convex, hence a closed interval.
Its endpoints are $z_-(u)$ and $z_+(u)$.
For $u,v\in\interior C$ and $0\le t\le1$, convexity of $D$ gives
\[
 \bigl((1-t)u+tv,(1-t)z_-(u)+tz_-(v)\bigr)\in D.
\]
Taking the minimum of $I_{(1-t)u+tv}$ yields
$z_-((1-t)u+tv)\le(1-t)z_-(u)+tz_-(v)$.
Taking upper endpoints instead gives the reverse inequality for $z_+$.
Thus both $z_-$ and $-z_+$ are convex. They are bounded because $D$ is compact.

Let $\varphi$ denote either of these two convex functions, and choose
$M_0>0$ such that $|\varphi|\le M_0$ on $\interior C$.
Fix $u\in\interior C$ and choose $\rho>0$ with
$\{v:\|v-u\|\le\rho\}\subset\interior C$.
For $0<\eta=\|v-u\|<\rho$, put $e=(v-u)/\eta$.
The identity $v=(1-\eta/\rho)u+(\eta/\rho)(u+\rho e)$ gives
\[
 \varphi(v)-\varphi(u)
 \le\frac{\eta}{\rho}\bigl(\varphi(u+\rho e)-\varphi(u)\bigr)
 \le\frac{2M_0\eta}{\rho}.
\]
Applying convexity to
$u=\frac{\rho}{\rho+\eta}v+\frac{\eta}{\rho+\eta}(u-\rho e)$ gives
\[
 (\rho+\eta)\varphi(u)
 \le\rho\varphi(v)+\eta\varphi(u-\rho e),
\]
and therefore
\[
 \varphi(v)-\varphi(u)
 \ge\frac{\eta}{\rho}\bigl(\varphi(u)-\varphi(u-\rho e)\bigr)
 \ge-\frac{2M_0\eta}{\rho}.
\]
These two estimates prove continuity.

To show that $z_-(u)<z_+(u)$ for $u\in\interior C$,
choose $d_*\in\interior D$ and $\varepsilon>0$ so that the open ball
of radius $\varepsilon$ about $d_*$ lies in $D$.
Write $u_*=\pi(d_*)$. Since $u'=(u-su_*)/(1-s)$ converges to $u\in\interior C$ as
$s\downarrow0$, choose $s\in(0,1)$ so small that $u'\in C$,
and then choose $d'\in D$ with $\pi(d')=u'$.
For $d=sd_*+(1-s)d'$ we have $\pi(d)=u$. If $\|b\|<s\varepsilon$, then
\[
 d+b=s\left(d_*+\frac{b}{s}\right)+(1-s)d'\in D.
\]
Writing $d=(u,z)$, choose $0<\delta<s\varepsilon$. Then
$(u,z-\delta),(u,z+\delta)\in D$, so
$z_-(u)\le z-\delta<z+\delta\le z_+(u)$.
Neither $(u,z_-(u))$ nor $(u,z_+(u))$ can lie in $\interior D$:
an open ball about either point would contain points of $D$ with
the same first coordinate and a last coordinate beyond the corresponding
minimum or maximum. Since $D$ is closed, both points belong to $\partial D$.

\smallskip
\noindent\textup{(2)} For $b>0$, part \textup{(1)} gives
\[
 \ell(u,z_-(u))
 =\ell(u,z_+(u))-b\bigl(z_+(u)-z_-(u)\bigr)<0.
\]
For $b<0$, it gives
\[
 \ell(u,z_+(u))
 =\ell(u,z_-(u))+b\bigl(z_+(u)-z_-(u)\bigr)<0.
\]
If $b=0$, then $a\ne0$ and $\langle a,u\rangle\le0$ on $C$.
Equality at an interior point $u$ is impossible: for sufficiently small
$t>0$, the point $u+ta\in C$ would satisfy
$\langle a,u+ta\rangle=t\|a\|^2>0$.
Hence $\ell(u,z_-(u))=\ell(u,z_+(u))=\langle a,u\rangle<0$.
\end{proof}

\subsection{Sections forced by proper extensions}

For $F:X\to\mathbb R^k$ and $U\subset\mathbb R^k$, a section
of $F$ on $U$ is a continuous map $s:U\to X$ satisfying
$F(s(u))=u$ for all $u\in U$. A proper admissible extension yields
such a section on $\interior C$ in the countably compact case and
on a radial interval in the locally compact case.

\begin{lemma}\label{lem:countable-evaluation}
Let $X$ be countably compact Hausdorff with at least two points,
and let $f_1,\ldots,f_k$ be a basis of an admissible subspace of $\C(X)$.
Put $F=(f_1,\ldots,f_k)$ and $C=\conv F(X)$.
Then $C$ is compact with nonempty interior, and, for
$Z=F^{-1}(\partial C)$, the map $F|_Z:Z\to\partial C$ is a homeomorphism.
If this subspace has a proper admissible extension in $\C(X)$,
there is a continuous $s:\interior C\to X$ satisfying $F(s(u))=u$.
\end{lemma}

\begin{proof}
A continuous image of a countably compact space is countably compact,
and a metrizable countably compact space is compact. Hence $F(X)$
is compact, and so is its convex hull by Lemma~\ref{lem:evaluation}.
If $C$ had empty interior, it would lie in an affine hyperplane:
$\sum_{j=1}^k a_jf_j=c$ for some $a\ne0$ and $c\in\mathbb R$.
For $c=0$ this contradicts linear independence; for $c\ne0$ it gives
a nonzero constant in an admissible subspace. Thus $\interior C\ne\varnothing$.

For each $a\ne0$, let $x_a$ be the unique maximizer of
$\langle a,F\rangle$. Lemma~\ref{lem:evaluation} gives
\[
 \{y\in C:\langle a,y\rangle=H_C(a)\}=\{F(x_a)\},
 \qquad F^{-1}(\{F(x_a)\})=\{x_a\}.
\]
Every boundary point of $C$ lies on a supporting hyperplane.
Consequently, $F|_Z$ is a continuous bijection onto $\partial C$.
If $H$ is closed in $Z$, then it is closed in $X$, because $Z$ is
closed in $X$. Thus $H$ is countably compact and its metrizable
image $F(H)$ is compact, hence closed in $\partial C$.
This makes $F|_Z$ a closed bijection, proving the homeomorphism.

Suppose that $g\notin\operatorname{span}\{f_1,\ldots,f_k\}$ gives
an admissible extension. Apply the preceding conclusions to
$\Psi=(F,g):X\to\mathbb R^{k+1}$ and $D=\conv\Psi(X)$.
Then $D$ is compact with nonempty interior and
$\Psi|_{Z_D}:Z_D\to\partial D$ is a homeomorphism, where
$Z_D=\Psi^{-1}(\partial D)$.
For $\pi(u,z)=u$, we have $\pi(D)=C$.
By Lemma~\ref{lem:vertical}, the map
$u\mapsto(u,z_+(u))$ is continuous from $\interior C$ into $\partial D$.
Therefore
\[
 s(u)=(\Psi|_{Z_D})^{-1}(u,z_+(u)),\qquad
 F(s(u))=\pi(\Psi(s(u)))=u
\]
defines the required continuous section.
\end{proof}

For functions in $\C_0(X)$, the sets
$\{x:\|F(x)\|\ge\delta\}$ are compact for every $\delta>0$.
This gives the compactness needed for the radial section.

\begin{lemma}\label{lem:radialsection}
Let $X$ be noncompact locally compact Hausdorff, and let
$E\subset\Chat_0(X)\cup\{0\}$ have basis $f_1,\ldots,f_k$,
where $k\in\mathbb N$. Let $L,F,C$ be as in Lemma~\ref{lem:basis}.
Choose $c\in\interior C$, taking $c=0$ when $0\in\interior C$.
If $E$ has a proper linear extension in $\Chat(X)\cup\{0\}$,
there are $y\in\partial(C-c)$ and a continuous map
$s:(0,1)\to X$ such that
\begin{equation}\label{eq:radialsection}
 F(s(t))=c+ty\qquad(0<t<1).
\end{equation}
\end{lemma}

\begin{proof}
The upper boundary of the extended evaluation body can be recovered
over interior points of $C$, except possibly along one segment.
A radial interval avoiding that segment gives the required section.
Choose $g\in\C(X)\setminus E$ such that
$W=E+\spanop\{g\}\subset\Chat(X)\cup\{0\}$.
Both $g$ and $-g$ attain finite maxima, so $g$ is bounded:
\[
 -\max_X(-g)\le g(x)\le\max_Xg\qquad(x\in X).
\]
Define
\[
 \Psi:X\to\mathbb R^{k+1},\qquad \Psi(x)=(F(x),g(x)),
 \qquad K=\overline{\Psi(X)},\qquad D=\conv K,
\]
where the closure is taken in $\mathbb R^{k+1}$.
The boundedness of $F$ and $g$ makes $K$ compact; hence $D$ is compact
by Lemma~\ref{lem:evaluation}\textup{(1)}.

\smallskip
\noindent\textup{(1)} We first identify the points belonging to
$\overline{\Psi(X\setminus H)}$ for every compact subset $H\subset X$.
Let $\mathcal K(X)$ be the collection of compact subsets of $X$, and put
\begin{equation}\label{eq:limitset}
 T=\bigcap_{H\in\mathcal K(X)}\overline{\Psi(X\setminus H)}.
\end{equation}
Every set in this intersection is a nonempty compact subset of $K$.
For finitely many compact sets $H_1,\ldots,H_r$, their union is compact
and is not all of $X$. Therefore
\[
 \varnothing\ne
 \overline{\Psi\left(X\setminus\bigcup_{j=1}^r H_j\right)}
 \subset\bigcap_{j=1}^r\overline{\Psi(X\setminus H_j)}.
\]
The finite intersection property in compact $K$ gives $T\ne\varnothing$.

For $\varepsilon>0$, the set
$H_\varepsilon=\bigcup_{j=1}^k\{x:|f_j(x)|\ge\varepsilon\}$ is compact.
Every point of $\overline{\Psi(X\setminus H_\varepsilon)}$ has its first
$k$ coordinates in $[-\varepsilon,\varepsilon]$.
Since this holds for every $\varepsilon>0$, \eqref{eq:limitset} implies
$T\subset\{0\}\times\mathbb R$.
If $z\in K\setminus T$, some compact $H\subset X$ satisfies
$z\notin\overline{\Psi(X\setminus H)}$. The decomposition
$\Psi(X)=\Psi(H)\cup\Psi(X\setminus H)$ and compactness of $\Psi(H)$ give
\[
 K=\Psi(H)\cup\overline{\Psi(X\setminus H)},
\]
so $z\in\Psi(H)$. We have proved
\begin{equation}\label{eq:limitdecomposition}
 T\subset\{0\}\times\mathbb R,\qquad K=\Psi(X)\cup T.
\end{equation}
This union need not be disjoint.
Define
\[
 \lambda=\max\{z\in\mathbb R:(0,z)\in T\},\qquad p_+=(0,\lambda)\in T.
\]
Thus $p_+$ is the point of $T$ with largest last coordinate;
$\lambda$ need not equal $\sup_Xg$.

The set $D$ has nonempty interior. Indeed, if it lay in an affine
hyperplane, there would be $(a,b)\ne0$ and $d\in\mathbb R$ with
\[
 \sum_{j=1}^k a_jf_j(x)+bg(x)=d\qquad(x\in X).
\]
For $d=0$, this contradicts independence of $f_1,\ldots,f_k,g$.
For $d\ne0$, it puts a nonzero constant in $W$, which is impossible.
For the projection $\pi(u,z)=u$, the set $\pi(K)$ is compact and hence
closed; since it contains $F(X)$, it contains $\overline{F(X)}$.
Conversely, if $z\in K$, choose $x_m\in X$ with $\Psi(x_m)\to z$;
this is possible because the closure is taken in the metrizable space
$\mathbb R^{k+1}$, regardless of whether $X$ is metrizable.
Then $F(x_m)=\pi(\Psi(x_m))\to\pi(z)$, so $\pi(z)\in\overline{F(X)}$.
Thus $\pi(K)=\overline{F(X)}$.
Also $F(L)$ is closed and contains $F(X)$, while continuity of $F$ and
density of $X$ in $L$ give $F(L)\subset\overline{F(X)}$.
Hence $F(L)=\overline{F(X)}$, and consequently
\begin{equation}\label{eq:projectionbody}
 \pi(D)=\conv\pi(K)=\conv F(L)=C.
\end{equation}

\smallskip
\noindent\textup{(2)} Fix a direction $v=(a,b)\in\mathbb R^k\times\mathbb R$
with $b>0$, and define
\[
 q_v:X\to\mathbb R,\qquad
 q_v(x)=\langle a,F(x)\rangle+bg(x),\qquad H_v=H_D(v).
\]
Here $D_v=\{d\in D:\langle v,d\rangle=H_v\}$,
as in Lemma~\ref{lem:evaluation}.
The coefficient $b>0$ and independence of the basis imply $q_v\ne0$.
Thus $q_v$ has exactly one maximizer $x_v\in X$.
The maximum of the linear functional $z\mapsto\langle v,z\rangle$
on $K=\overline{\Psi(X)}$ equals $\sup_Xq_v=\max_Xq_v$.
It also equals its maximum on $D=\conv K$. Since $p_+\in K$, we obtain
\[
 H_v=\max_Xq_v\ge\langle v,p_+\rangle=b\lambda.
\]
If $H_v>b\lambda$, no point of $T$ is a maximizer, by
\eqref{eq:limitdecomposition}. The only maximizing point in $K$ is
$\Psi(x_v)$.
If $H_v=b\lambda$, a point $(0,z)\in T$ is maximizing precisely when
$bz=b\lambda$, hence precisely when it is $p_+$.
The only maximizing evaluation from $X$ is still $\Psi(x_v)$.
Applying \eqref{eq:evaluationface} to the inclusion of $K$ into
$\mathbb R^{k+1}$ gives
\begin{equation}\label{eq:upperfaces}
 D_v=
 \begin{cases}
 \{\Psi(x_v)\},&H_v>b\lambda,\\
 [p_+,\Psi(x_v)],&H_v=b\lambda.
 \end{cases}
\end{equation}

\smallskip
\noindent\textup{(3)} All maximizing points $x_v$ arising in the equality
case of \eqref{eq:upperfaces} coincide.
To prove this, let $v=(a,b)$ and $v'=(a',b')$ have $b,b'>0$ and
$H_v=b\lambda$, $H_{v'}=b'\lambda$.
The functions
$q_v-b\lambda$ and $q_{v'}-b'\lambda$ are nonpositive on $X$.
Their sum has supremum zero, since the corresponding affine functional
is nonpositive on $K$ and vanishes at $p_+$:
\[
 \begin{aligned}
 \sup_{x\in X}\bigl(q_v(x)+q_{v'}(x)-(b+b')\lambda\bigr)
 &=\max_{z\in K}\bigl(\langle v+v',z\rangle-(b+b')\lambda\bigr)\\
 &=0.
 \end{aligned}
\]
Since $b+b'>0$, the function $q_v+q_{v'}$ is a nonzero member of $W$.
It therefore attains its maximum $(b+b')\lambda$ at a point $x\in X$.
At that point the two nonpositive summands have sum zero, so
\[
 q_v(x)-b\lambda=0,\qquad q_{v'}(x)-b'\lambda=0.
\]
Uniqueness of each maximizer gives $x=x_v=x_{v'}$.
If there exists $v=(a,b)$ with $b>0$ and $H_v=b\lambda$,
denote the common maximizer by $x_+$ and put $B_+=[0,F(x_+)]$.
If no such $v$ exists, put $B_+=\varnothing$.

\smallskip
\noindent\textup{(4)} By Lemma~\ref{lem:vertical} and
\eqref{eq:projectionbody}, the map
\[
 \Gamma:\interior C\to\partial D,\qquad
 \Gamma(u)=(u,z_+(u))
\]
is continuous. Put $\Omega=\interior C\setminus(B_+\cup\{0\})$.
We construct a continuous inverse to $\Psi$ along $\Gamma(\Omega)$.

Fix $u\in\interior C$. By the supporting hyperplane theorem,
choose $(a,b)\ne0$ such that $(w,z)\mapsto\langle a,w\rangle+bz$
attains its maximum on $D$ at $\Gamma(u)$. If $b<0$, the point $(u,z_-(u))\in D$ gives
\[
 \langle a,u\rangle+bz_-(u)
 -\bigl(\langle a,u\rangle+bz_+(u)\bigr)
 =-b\bigl(z_+(u)-z_-(u)\bigr)>0,
\]
contradicting the fact that $(u,z_+(u))$ maximizes
$(w,z)\mapsto\langle a,w\rangle+bz$ on $D$.
If $b=0$, then $a\ne0$ and $u$ would maximize
$w\mapsto\langle a,w\rangle$ on $C$. This is impossible at an interior
point, since $u+ta\in C$ for small $t>0$ and increases that value by
$t\|a\|^2$. Thus $b>0$.

For $u\in\Omega$, set $v=(a,b)$ for the direction just chosen.
The equality $H_v=b\lambda$ is impossible: by \eqref{eq:upperfaces},
it would imply $D_v=[p_+,\Psi(x_+)]$ and hence
$u=\pi(\Gamma(u))\in\pi(D_v)=B_+$, a contradiction.
Thus $H_v>b\lambda$ and $D_v=\{\Psi(x_v)\}$ by
\eqref{eq:upperfaces}. Since $\Gamma(u)\in D_v$,
$\Gamma(u)=\Psi(x_v)$. This point has exactly one preimage
in $X$, because every preimage would maximize $q_v$.
Define
\begin{equation}\label{eq:localsigma}
 \sigma:\Omega\to X,\qquad
 \Psi(\sigma(u))=\Gamma(u),\qquad F(\sigma(u))=u.
\end{equation}

To prove continuity, let $Q\subset\Omega$ be a closed Euclidean ball.
Since $0\notin Q$, put $\delta=\min_{u\in Q}\|u\|>0$.
The set $\Psi^{-1}(\Gamma(Q))$ is closed in $X$ and is contained in
\[
 \{x\in X:\|F(x)\|\ge\delta\}
 \subset\bigcup_{j=1}^k
       \{x\in X:|f_j(x)|\ge\delta/\sqrt{k}\}.
\]
The last inclusion follows because $|f_j(x)|<\delta/\sqrt{k}$ for
all $j$ would give $\|F(x)\|^2<k(\delta^2/k)=\delta^2$.
The union on the right is compact, since each $f_j$ vanishes at infinity.
Hence $\Psi^{-1}(\Gamma(Q))$ is compact.
The restriction of $\Psi$ to this set is a continuous bijection onto
$\Gamma(Q)$, and so has continuous inverse.
Composing that inverse with $\Gamma|_Q$ proves continuity of $\sigma|_Q$.
The set $\Omega$ is open; such balls are neighborhoods locally at every
point of $\Omega$. Thus $\sigma$ is continuous.

\smallskip
\noindent\textup{(5)} We now choose $y\in\partial(C-c)$ such that
$\{c+ty:0<t<1\}\subset\Omega$.
Let $p:\mathbb R^k\to[0,\infty)$ be the Minkowski functional of $C-c$ from
Lemma~\ref{lem:gauge}.
Suppose first that $0\in\interior C$, so $c=0$.
If $B_+$ is nondegenerate, write $v_+=F(x_+)\ne0$ and choose
$y=-v_+/p(-v_+)\in\partial C$.
For $0<t<1$, the point $ty$ is nonzero and lies on the ray opposite to
$B_+=[0,v_+]$. Also $p(ty)=t<1$, so $ty\in\interior C$.
Hence $ty\in\Omega$.
If $B_+$ is empty or equals $\{0\}$, take any $y\in\partial C$;
the same conclusion follows from $p(ty)=t$ and $ty\ne0$.
The map $s(t)=\sigma(ty)$ is continuous and satisfies
\eqref{eq:radialsection} by \eqref{eq:localsigma}.

Finally, suppose $0\in\partial C$.
By the supporting hyperplane theorem, choose $a_0\ne0$ such that
$\langle a_0,u\rangle\le0$ for every $u\in C$, and put
$f_0=\sum_{j=1}^k(a_0)_jf_j$.
Lemma~\ref{lem:basis}\textup{(3)} gives
$f_0\le0$ and a unique zero $x_0\in X$.
For a vector $v=(a,b)$ with $b>0$ and $H_v=b\lambda$, both $f_0$ and
$q_v-b\lambda$ are nonpositive. Their sum has supremum zero: its
associated affine functional on $K$ vanishes at $p_+$.
The function $f_0+q_v$ is a nonzero member of $W$, because its coefficient
of $g$ is $b>0$. It attains the maximum $b\lambda$ at some $x\in X$.
Then
\[
 0=f_0(x)+\bigl(q_v(x)-b\lambda\bigr),\qquad
 f_0(x)\le0,\quad q_v(x)-b\lambda\le0,
\]
so $f_0(x)=0$ and $q_v(x)=b\lambda$. Thus $x=x_0=x_v$.
Consequently, if such a vector $v$ exists, $x_+=x_0$ and
\[
 B_+=[0,F(x_0)]\subset
 \{u\in C:\langle a_0,u\rangle=0\}\subset\partial C.
\]
If none exists, $B_+=\varnothing$. In either case,
$\Omega=\interior C$ because $0\in\partial C$.
Choose any $y\in\partial(C-c)$. Lemma~\ref{lem:gauge} gives
$p(c+ty-c)=tp(y)=t<1$, so $c+ty\in\Omega$ for $0<t<1$.
Now $s(t)=\sigma(c+ty)$ is the required map.
\end{proof}

\subsection{A perturbation preserving all extrema}

A continuous scalar function, arbitrarily close to the identity
and fixing the endpoints, can fail to admit a continuous right inverse
on $(0,1)$.

\begin{lemma}\label{lem:fold}
The following statements hold.
\begin{enumerate}
\item For every $\varepsilon>0$ there is a continuous function
$w:[0,1]\to[0,1]$ such that
\begin{equation}\label{eq:foldconditions}
 w(0)=0,\qquad w(1)=1,\qquad 0<w(t)<1\quad(0<t<1),
\end{equation}
$\|w-\mathrm{id}_{[0,1]}\|_\infty<\varepsilon$, and $w$ is not strictly
increasing.
\item Every continuous $w:[0,1]\to[0,1]$ satisfying
\eqref{eq:foldconditions} and admitting a continuous function
$r:(0,1)\to[0,1]$ with $w(r(t))=t$ for $0<t<1$ is strictly increasing
on $[0,1]$.
\end{enumerate}
\end{lemma}

\begin{proof}
\textup{(1)} Choose $0<\tau<\min\{\varepsilon,1/6\}$ and define
\[
 w(t)=
 \begin{cases}
 t,&0\le t\le\frac12,\\
 2t-\frac12,&\frac12\le t\le\frac12+\tau,\\
 1+3\tau-t,&\frac12+\tau\le t\le\frac12+2\tau,\\
 2t-\frac12-3\tau,&\frac12+2\tau\le t\le\frac12+3\tau,\\
 t,&\frac12+3\tau\le t\le1.
 \end{cases}
\]
The adjacent formulas agree at each endpoint. Their values at the four
interior breakpoints are, in order,
$1/2$, $1/2+2\tau$, $1/2+\tau$, and $1/2+3\tau$.
They lie in $(0,1)$, so piecewise affinity proves
\eqref{eq:foldconditions}.
On the middle interval, the slope is $-1$, hence $w$ is not strictly
increasing. On each interval, $w(t)-t$ is affine; its values at all
breakpoints, including $0$ and $1$, are
$0,0,\tau,-\tau,0,0$. Therefore
$\|w-\mathrm{id}_{[0,1]}\|_\infty=\tau<\varepsilon$.

\smallskip
\noindent\textup{(2)} The identity $w(r(t))=t$ makes $r$ injective, and
\eqref{eq:foldconditions} gives $0<r(t)<1$.
For $0<\delta<1$, compactness and \eqref{eq:foldconditions} imply
\[
 m_\delta=\min_{\delta\le u\le1}w(u)>0,\qquad
 M_\delta=\max_{0\le u\le1-\delta}w(u)<1.
\]
If $0<t<m_\delta$, then $r(t)<\delta$, since otherwise
$t=w(r(t))\ge m_\delta$.
If $M_\delta<t<1$, then $r(t)>1-\delta$, since otherwise
$t=w(r(t))\le M_\delta$.
Thus
\[
 \lim_{t\downarrow0}r(t)=0,\qquad
 \lim_{t\uparrow1}r(t)=1.
\]
A continuous injective real function on an interval is strictly monotone.
The two limits force $r$ to be increasing, and the intermediate value
theorem gives $r((0,1))=(0,1)$.
For $0<a<b<1$, write $a=r(t_a)$ and $b=r(t_b)$.
Then $t_a<t_b$, whence
$w(a)=t_a<t_b=w(b)$.
The endpoint comparisons $w(0)<w(b)$ and $w(a)<w(1)$ follow from
\eqref{eq:foldconditions}. Thus $w$ is strictly increasing on $[0,1]$.
\end{proof}

The deformation of the evaluation map yields an admissible
isometric image of the original subspace with no proper admissible
extension. This image need not contain the original subspace, since
the definition of $(\alpha,\beta)$-lineability concerns every
prescribed $\alpha$-dimensional subspace.

\begin{theorem}\label{lem:perturbation}
Let $X$ be Hausdorff with at least two points, and let $E$ be a
nonzero finite-dimensional admissible subspace. Suppose either that
$X$ is locally compact and $E\subset\C_0(X)$, or that $X$ is
countably compact and $E\subset\C(X)$.
For every $\eta>0$ there is a linear map $S:E\to\C_b(X)$ with the
following properties. In the locally compact case with
$E\subset\C_0(X)$, one may also require $S(E)\subset\C_0(X)$.
\begin{enumerate}
\item For every $f\in E\setminus\{0\}$,
\[
 \begin{aligned}
 \max_X Sf&=\max_Xf,&\qquad \min_X Sf&=\min_Xf,\\
 \Argmax_X Sf&=\Argmax_Xf,&
 \Argmin_X Sf&=\Argmin_Xf.
 \end{aligned}
\]
\item The map $S$ is an isometry and $\|S-I_E\|<\eta$, where
$I_E:E\to\C_b(X)$ is the inclusion and $E$ carries the supremum norm.
\item Every linear subspace $W$ satisfying
$S(E)\subset W\subset\Chat(X)\cup\{0\}$ equals $S(E)$.
\end{enumerate}
\end{theorem}

\begin{proof}
\textup{(1)} Choose a basis $f_1,\ldots,f_k$ of $E$.
If $X$ is locally compact and $E\subset\C_0(X)$, use the setup
$L,F,C$ of Lemma~\ref{lem:basis}; call this case (i).
Otherwise use the countably compact hypothesis and set $L=X$,
$F=(f_1,\ldots,f_k)$ and $C=\conv F(X)$; call this case (ii).
Lemma~\ref{lem:countable-evaluation} applies in case (ii).
In both cases $F(L)$ and $C$ are compact and $\interior C\ne\varnothing$.
Choose $c\in\interior C$, taking $c=0$ in case (i) when $X$ is
noncompact and $0\in\interior C$.
Let $p:\mathbb R^k\to[0,\infty)$ be the Minkowski functional of $C-c$ from
Lemma~\ref{lem:gauge}. Choose $w$ by Lemma~\ref{lem:fold}\textup{(1)} so that
\[
 \varepsilon:=\|w-\mathrm{id}_{[0,1]}\|_\infty<\eta/2.
\]
Define
\begin{equation}\label{eq:gaugedeformation}
 \begin{aligned}
 h:L\to[0,1],&\qquad h(x)=p(F(x)-c),\\
 G:L\to\mathbb R^k,&\qquad
 G(x)=
 \begin{cases}
 c+\dfrac{w(h(x))}{h(x)}(F(x)-c),&h(x)>0,\\[5pt]
 c,&h(x)=0.
 \end{cases}
 \end{aligned}
\end{equation}
The range of $h$ lies in $[0,1]$ by \eqref{eq:gaugesets}.
Continuity of $F$ and $p$ gives continuity of $h$.
When $h(x)>0$, positive homogeneity gives
\[
 p\left(\frac{F(x)-c}{h(x)}\right)=1,
\]
so $(F(x)-c)/h(x)\in\partial(C-c)$.
Choose $R_0>0$ with $\|y\|\le R_0$ on $C-c$.
At every point with $h(x)>0$, we have
\[
 \|G(x)-c\|
 =w(h(x))\left\|\frac{F(x)-c}{h(x)}\right\|
 \le R_0w(h(x)).
\]
Since $w(0)=0$, this bound proves continuity of $G$ also where $h=0$.
For $h>0$, homogeneity of $p$ gives
$p(G(x)-c)=w(h(x))p((F(x)-c)/h(x))=w(h(x))$.
For $h=0$, both sides are zero. Thus
\begin{equation}\label{eq:gaugeprofile}
 p(G(x)-c)=w(h(x))\qquad(x\in L).
\end{equation}
In particular, $G(L)\subset C$.
If $F(x)\in\partial C$, then $h(x)=1$ and $w(h(x))=1$, hence $G(x)=F(x)$.
If $F(x)\in\interior C$, then $h(x)<1$ and $w(h(x))<1$, hence
$G(x)\in\interior C$. We have proved
\begin{equation}\label{eq:boundaryfixed}
 \begin{aligned}
 F(x)\in\partial C&\ \Longrightarrow\ G(x)=F(x),\\
 F(x)\in\interior C&\ \Longrightarrow\ G(x)\in\interior C.
 \end{aligned}
\end{equation}

Write $G=(g_1,\ldots,g_k)$ and restrict the coordinates to $X$.
These restrictions are continuous and bounded, since $G(L)\subset C$.
In case (ii) this is the required conclusion. In case (i) we also
verify that they belong to $\C_0(X)$. This is automatic for compact $X$.
For noncompact $X$ in case (i), $F(\infty)=0\in C$.
If $0\in\interior C$, our choice $c=0$ gives $h(\infty)=0$ and
$G(\infty)=0$.
If $0\in\partial C$, the first implication in \eqref{eq:boundaryfixed}
gives $G(\infty)=F(\infty)=0$.
For each $j$ and $\delta>0$, the set
$\{x\in L:|g_j(x)|\ge\delta\}$ is therefore compact in $L$ and avoids
$\infty$. It is a compact subset of $X$, proving $g_j|_X\in\C_0(X)$.
We use the same symbol $g_j$ for the restriction.

Define the linear map
\[
 S:E\to\C_b(X),\qquad
 S\left(\sum_{j=1}^k a_jf_j\right)=\sum_{j=1}^k a_jg_j.
\]
This is well defined because the representation in the basis is unique.
Fix $a\in\mathbb R^k\setminus\{0\}$ and put
$f=\sum_{j=1}^k a_jf_j$.
Lemma~\ref{lem:evaluation} gives
$H_C(a)=\max_Xf$.
An interior point of $C$ cannot maximize the nonzero functional
$u\mapsto\langle a,u\rangle$, because a small displacement in the
direction $a$ increases its value while remaining in $C$.
On the boundary, $G$ agrees with $F$ by \eqref{eq:boundaryfixed}.
Consequently, for every $x\in X$,
\begin{equation}\label{eq:preservemaxset}
 \langle a,G(x)\rangle=H_C(a)
 \quad\Longleftrightarrow\quad
 \langle a,F(x)\rangle=H_C(a).
\end{equation}
Moreover, $\langle a,G(x)\rangle\le H_C(a)$ since $G(x)\in C$,
and equality is attained at the original maximizer of $f$.
Thus $\max_XSf=\max_Xf$ and $\Argmax_XSf=\Argmax_Xf$.
Applying \eqref{eq:preservemaxset} to $-a$ yields
$\max_X(-Sf)=\max_X(-f)$ and equality of their maximizing sets.
Since $\min f=-\max(-f)$ and $\Argmin f=\Argmax(-f)$, this proves all
four equalities in part \textup{(1)}.

\smallskip
\noindent\textup{(2)} For $f\ne0$, part \textup{(1)} gives
\[
 \|Sf\|_\infty
 =\max\{|\max_XSf|,|\min_XSf|\}
 =\max\{|\max_Xf|,|\min_Xf|\}
 =\|f\|_\infty.
\]
The identity also holds for $f=0$, so $S$ is an isometry.
In particular, it is injective and $S(E)\subset\Chat(X)\cup\{0\}$.

We prove the approximation estimate for
$f=\sum_{j=1}^k a_jf_j$.
For every $u\in C$, write $u=\sum_{i=1}^r\lambda_iF(x_i)$ with
$\lambda_i\ge0$ and $\sum_{i=1}^r\lambda_i=1$. Then
\begin{equation}\label{eq:functionalbound}
 |\langle a,u\rangle|
 =\left|\sum_{i=1}^r\lambda_if(x_i)\right|
 \le\sum_{i=1}^r\lambda_i\|f\|_\infty
 =\|f\|_\infty.
\end{equation}
In case (i), when $X$ is noncompact, this also holds at $x_i=\infty$, where $f=0$.
If $h(x)>0$, put $y=(F(x)-c)/h(x)$. Then $c+y\in\partial C\subset C$
and $c\in C$, so \eqref{eq:functionalbound} gives
\[
 |\langle a,y\rangle|
 \le|\langle a,c+y\rangle|+|\langle a,c\rangle|
 \le2\|f\|_\infty.
\]
Subtracting $F(x)=c+h(x)y$ from $G(x)=c+w(h(x))y$ now gives
\[
 \begin{aligned}
 |Sf(x)-f(x)|
 &=|\langle a,G(x)-F(x)\rangle|\\
 &=|w(h(x))-h(x)|\,|\langle a,y\rangle|\\
 &\le2\varepsilon\|f\|_\infty.
 \end{aligned}
\]
If $h(x)=0$, then $F(x)=G(x)=c$ and the difference is zero.
Taking first the supremum over $x\in X$ and then over
$\|f\|_\infty\le1$ gives
\[
 \|S-I_E\|
 =\sup_{\|f\|_\infty\le1}\sup_{x\in X}|Sf(x)-f(x)|
 \le2\varepsilon<\eta.
\]

\smallskip
\noindent\textup{(3)} We first show that
\begin{equation}\label{eq:samebody}
 \conv G(L)=C.
\end{equation}
The set $G(L)$ is compact: in case (i), $L$ is compact; in case (ii),
$L$ is countably compact and $G$ has finite-dimensional range.
Thus $\conv G(L)$ is compact. Its inclusion in $C$ follows from $G(L)\subset C$.
For each $a\ne0$, the maximum of $f=\langle a,F\rangle$ on $X$ equals
$H_C(a)$ and is attained at a point fixed by the deformation.
Hence
\[
 \max_{u\in\conv G(L)}\langle a,u\rangle
 =\max_{x\in L}\langle a,G(x)\rangle=H_C(a).
\]
If $u_0\in C\setminus\conv G(L)$, separation of the point $u_0$ from the
compact convex set $\conv G(L)$ would give $a\ne0$ with
\[
 H_C(a)\ge\langle a,u_0\rangle
 >\max_{u\in\conv G(L)}\langle a,u\rangle=H_C(a),
\]
a contradiction. This proves \eqref{eq:samebody}.

Suppose now that $S(E)$ has a proper extension in
$\Chat(X)\cup\{0\}$.
In case (i) with noncompact $X$, apply Lemma~\ref{lem:radialsection} to the basis
$g_1,\ldots,g_k$ of $S(E)$, using \eqref{eq:samebody} and the chosen $c$.
This gives $y\in\partial(C-c)$ and a continuous $s:(0,1)\to X$ with
$G(s(t))=c+ty$.
In all other cases $X$ is countably compact, so
Lemma~\ref{lem:countable-evaluation}, applied to the basis
$g_1,\ldots,g_k$, gives a continuous section on $\interior C$.
Composing it with $t\mapsto c+ty$, for any $y\in\partial(C-c)$,
gives the same conclusion.
The construction of $G$ replaces $h(x)=p(F(x)-c)$ by $w(h(x))$.
The map $s$ just obtained would yield a continuous function $r$
satisfying $w(r(t))=t$ for $0<t<1$. Indeed, define $r:(0,1)\to[0,1]$ by $r(t)=h(s(t))$.
This map is continuous, and \eqref{eq:gaugeprofile} gives
\[
 w(r(t))
 =p(G(s(t))-c)
 =p(ty)
 =tp(y)
 =t\qquad(0<t<1).
\]
Lemma~\ref{lem:fold}\textup{(2)} would make $w$ strictly increasing,
contrary to its construction in part \textup{(1)} of that lemma.
Thus no proper extension exists, proving part \textup{(3)}.
\end{proof}

\begin{proof}[Proof of Theorem~\ref{thm:A}]
Let $1\le\alpha<\aleph_0$ and $\beta\ge\alpha$.
If $\Chat_0(X)$ is not $\alpha$-lineable, both $(\alpha,\beta)$-lineability and $(\alpha,\beta)$-spaceability fail by their definitions.
Suppose it is $\alpha$-lineable, and choose
$E\subset\Chat_0(X)\cup\{0\}$ with $\dim E=\alpha$.
Theorem~\ref{lem:perturbation} gives an isometric image $S(E)$ of dimension
$\alpha$ which has no proper admissible linear extension.
If $\beta>\alpha$, any extension of dimension $\beta$ would be proper.
Therefore neither $(\alpha,\beta)$-lineability nor
$(\alpha,\beta)$-spaceability holds.
If $\beta=\alpha$, every prescribed $\alpha$-dimensional subspace serves
as its own extension. It is closed because finite-dimensional subspaces
of normed spaces are closed.
For countably compact Hausdorff $X$, apply the countably compact
case of Theorem~\ref{lem:perturbation} to an admissible subspace of
$\C(X)$. The same argument excludes every $\beta>\alpha$.
On the diagonal, finite-dimensional subspaces are closed in both
ambient topologies, since both are Hausdorff topological vector topologies.
\end{proof}

\begin{corollary}\label{cor:transfer}
Let $X$ be locally compact Hausdorff with at least two points, and
let $1\le\alpha<\aleph_0$. Assume that $\Chat_0(X)$ is
$\alpha$-lineable. For every cardinal $\beta\ge\alpha$,
\[
 \Chat(X)\text{ is }(\alpha,\beta)\text{-lineable}
 \quad\Longleftrightarrow\quad \alpha=\beta.
\]
The same equivalence holds for $(\alpha,\beta)$-spaceability in
$\C(X)$ with the compact-open topology and in $\C_b(X)$ with the
supremum norm.
\end{corollary}

\begin{proof}
Choose $E\subset\Chat_0(X)\cup\{0\}$ with $\dim E=\alpha$.
Theorem~\ref{lem:perturbation} gives an isometric image $S(E)$ admitting
no proper admissible extension in $\C(X)$. This excludes
$(\alpha,\beta)$-lineability, and hence both spaceability properties,
when $\beta>\alpha$.
For $\beta=\alpha$, every prescribed subspace is its own extension
and is closed in both ambient Hausdorff topological vector spaces,
since its dimension is finite. The required $\alpha$-lineability
follows from $E\subset\Chat(X)\cup\{0\}$.
\end{proof}

\section{Proof of Theorem B}\label{sec:proofB}

By Theorem~\ref{thm:A}, it remains to determine the admissible
dimensions on $M$. For compact $M$, these are characterized by embedded
spheres; the noncompact case requires a separate argument at infinity.

\subsection{Sphere zero-sets and finite-dimensional subspaces}

The criterion combines the sphere construction of
\cite[Theorem 2.1]{BCFPS} with the evaluation-map method of
\cite[Lemmas 3.2--3.5]{BCFPS}.

\Needspace{16\baselineskip}
\begin{proposition}\label{prop:countably-compact-sphere}
Let $X$ be a countably compact space with at least two points whose
continuous real functions separate points, and let $k\geq1$ be an integer. The following assertions
are equivalent:
\begin{enumerate}
    \item $\Chat(X)$ is $k$-lineable;
    \item there exists a continuous function $u:X\to[0,1]$ such that
    $u^{-1}(0)$ is homeomorphic to $S^{k-1}$.
\end{enumerate}
Whenever these conditions hold, there is a linear isometry
$T:\ell_2^k\to \C_b(X)$ such that
$T(\ell_2^k)\setminus\{0\}\subset\Chat(X)$,
where $\ell_2^k$ denotes $\mathbb R^k$ with its Euclidean norm.
Moreover, the implication {\rm (2)}$\Rightarrow${\rm (1)} holds
whenever continuous real functions on $X$ separate points, without
countable compactness.
\end{proposition}

\begin{proof}
Assume (1), choose a basis $f_1,\ldots,f_k$ of an admissible subspace,
and put $F=(f_1,\ldots,f_k)$ and $C=\conv F(X)$.
The separation assumption implies that $X$ is Hausdorff.
Lemma~\ref{lem:countable-evaluation} gives a homeomorphism from
$Z=F^{-1}(\partial C)$ onto $\partial C$, and
Lemma~\ref{lem:gauge} gives $\partial C\cong S^{k-1}$.
The continuous function
\[
 u(x)=\min\{1,\operatorname{dist}(F(x),\partial C)\}
\]
has $u^{-1}(0)=Z$, proving (2).

Conversely, assume (2) and that continuous real functions separate
points of $X$. Put $Z=u^{-1}(0)$ and choose a homeomorphism
$\phi:Z\to S^{k-1}$. The set $Z$ is compact.
We extend the coordinate functions of $\phi$ through an
evaluation map into a cube.
Let $I$ be the set of continuous functions $f:X\to[0,1]$,
and define the evaluation map
\[
 e:X\longrightarrow[0,1]^I,\qquad e(x)=(f(x))_{f\in I}.
\]
It is injective: a continuous real function distinguishing two points
can be composed with the continuous injective function
$t\mapsto\tfrac12+\tfrac1\pi\arctan t$ to give a member of $I$
distinguishing them. Only the restriction of $e$ to $Z$ needs to be a topological
embedding. Indeed, $e|_Z$ is a continuous bijection
from a compact space onto a Hausdorff space, hence a homeomorphism.
The set $e(Z)$ is compact and closed in the compact Hausdorff cube
$[0,1]^I$. By the Tietze extension theorem on this cube, each
$\phi_j\circ(e|_Z)^{-1}$ extends to a continuous
$\widetilde A_j:[0,1]^I\to[-1,1]$.
Then $A_j=\widetilde A_j\circ e$ extends $\phi_j$ to $X$.
Writing $A=(A_1,\ldots,A_k):X\to\mathbb R^k$, define
\[
 G:X\to\mathbb R^k,\qquad
 G(x)=\frac{(1-u(x))A(x)}{\max\{1,\|A(x)\|\}},
\]
and define $T:\mathbb R^k\to\C_b(X)$ by
$(Ta)(x)=\langle a,G(x)\rangle$ for $a\in\mathbb R^k$ and
$x\in X$.
The map $G$ is continuous, $G|_Z=\phi$, and
$\|G(x)\|\leq1-u(x)\leq1$ for every $x\in X$.
Thus $T:\ell_2^k\to \C_b(X)$ is linear and
$\|Ta\|_\infty\leq\|a\|$.

Fix $a\ne0$ and let $z_a=\phi^{-1}(a/\|a\|)$.
Then $(Ta)(z_a)=\|a\|$.
For $x\notin Z$, we have $u(x)>0$ and therefore
\[
(Ta)(x)\leq|\langle a,G(x)\rangle|
\leq\|a\|\|G(x)\|
\leq\|a\|(1-u(x))<\|a\|.
\]
For $x\in Z$, both $\phi(x)$ and $a/\|a\|$ are unit vectors,
and
\[
\left\|\phi(x)-\frac{a}{\|a\|}\right\|^2
=2-\frac{2\langle a,\phi(x)\rangle}{\|a\|}.
\]
Hence $(Ta)(x)=\|a\|$ holds precisely when
$\phi(x)=a/\|a\|$, equivalently when $x=z_a$.
Thus $Ta$ has a unique point of maximum and
$\|Ta\|_\infty=\|a\|$.
This proves both (1) and the asserted isometric conclusion.
No countable compactness was used in this implication.
\end{proof}

On a compact metrizable space, every embedded sphere is a zero-set,
so Proposition~\ref{prop:countably-compact-sphere} simplifies as follows.

\begin{corollary}\label{cor:compact-metric-sphere}
Let $K$ be a compact metrizable space with at least two points,
and let $k\geq1$ be an integer. Then $\Chat(K)$ is $k$-lineable
if and only if $K$ contains a subset homeomorphic to $S^{k-1}$.
Whenever these conditions hold, there is a linear isometry
$T:\ell_2^k\to\C(K)$ such that
$T(\ell_2^k)\setminus\{0\}\subset\Chat(K)$.
\end{corollary}

\begin{proof}
If $\Chat(K)$ is $k$-lineable,
Proposition~\ref{prop:countably-compact-sphere} gives a zero-set
in $K$ homeomorphic to $S^{k-1}$, which is in particular a
subset of $K$ with the required property.
Conversely, suppose that $Z\subset K$ is homeomorphic to
$S^{k-1}$. Since the sphere is compact, $Z$ is compact and
therefore closed in $K$. Choose a compatible metric $d$ on $K$
and define
\[
u(x)=\min\{1,d(x,Z)\},\qquad
 d(x,Z)=\inf_{z\in Z}d(x,z).
\]
The inequality $|d(x,Z)-d(y,Z)|\leq d(x,y)$ shows that $u$
is continuous. Since $Z$ is closed, $u(x)=0$ holds precisely
when $x\in Z$. Thus $u^{-1}(0)=Z$, and the proposition gives
both $k$-lineability and the asserted isometry, with
$\C_b(K)=\C(K)$ by compactness.
\end{proof}

\begin{example}\label{ex:noncompact-countably-compact}
Let $\omega_1$ be the first uncountable ordinal and equip
$L=[0,\omega_1)$ with the order topology.

The ordinal space $L$ is normal Hausdorff and countably compact,
but not compact. Indeed, every countable subset is contained in a compact
ordinal interval $[0,\beta]$ with $\beta<\omega_1$, whereas the open cover
$\{[0,\alpha):0<\alpha<\omega_1\}$ has no finite subcover.

For each integer $k\geq1$, consider the topological
disjoint union
\[
X_k=S^{k-1}\sqcup L.
\]
Each component is open and closed in $X_k$.
Since a finite topological disjoint union of normal
Hausdorff countably compact spaces has the same
properties, $X_k$ satisfies the hypotheses of
Proposition~\ref{prop:countably-compact-sphere}.
However, $X_k$ is not compact, because its closed
subspace $L$ is not compact.

The function $u:X_k\to[0,1]$ defined by $u=0$ on
$S^{k-1}$ and $u=1$ on $L$ is continuous and satisfies
$u^{-1}(0)=S^{k-1}$. Hence the proposition yields
$k$-lineability of $\Chat(X_k)$.

In this example the corresponding linear isometry
can also be written explicitly:
\[
T:\ell_2^k\longrightarrow \C_b(X_k),\qquad
(Ta)(x)=
\begin{cases}
\langle a,x\rangle,&x\in S^{k-1},\\
0,&x\in L.
\end{cases}
\]
The function $Ta$ is continuous because both components
are open, and $\|Ta\|_\infty=\|a\|$.
For $a\ne0$, its maximum is $\|a\|>0$ and is attained
only at $a/\|a\|\in S^{k-1}$. Indeed, on the sphere,
\[
\|a\|-\langle a,x\rangle
=\frac{\|a\|}{2}
 \left\|x-\frac{a}{\|a\|}\right\|^2,
\]
whereas $Ta=0$ on $L$.
Thus $T(\ell_2^k)\setminus\{0\}\subset\Chat(X_k)$.

\end{example}

\begin{remark}\label{rem:compactness-zero-set}
The countable compactness assumption in
Proposition~\ref{prop:countably-compact-sphere}
cannot simply be omitted: the implication
\textup{(1)$\Rightarrow$(2)} fails even for
$X=\mathbb R$ and $k=2$.
By \cite[Theorem 16]{GQ}, $\Chat_0(\mathbb R)$ is
$2$-lineable, so $\Chat(\mathbb R)$ is also $2$-lineable.
However, $\mathbb R$ contains no subset homeomorphic to $S^1$.
Indeed, such a subset would be compact and connected,
hence a nondegenerate closed interval. Removing an interior
point disconnects that interval, whereas removing any
single point from $S^1$ leaves a connected space.
Therefore $\mathbb R$ contains no zero-set homeomorphic
to $S^1$, and the equivalence in
Proposition~\ref{prop:countably-compact-sphere} does not extend to all normal
Hausdorff spaces.
\end{remark}

\subsection{Local lower bounds and countably compact manifolds}

A coordinate chart gives the lower bound, with all functions
supported in a common compact set.

\begin{lemma}\label{lem:chart}
Let $X$ be Hausdorff and suppose that some nonempty open subset of $X$
is homeomorphic to an open subset of $\mathbb R^n$, where $n\ge1$.
There exist a compact set $H\subset X$ and a linear isometry
$T:\ell_2^n\to\C_0(X)$ such that every $Ta$ vanishes outside $H$
and every $Ta$ with $a\ne0$ has a unique point of maximum.
In particular, $\Chat_0(X)$ is $n$-lineable.
\end{lemma}

\begin{proof}
Choose a homeomorphism $\chi:U\to O$, where $U\subset X$ and
$O\subset\mathbb R^n$ are nonempty open sets. Choose $a\in O$ and $r>0$ such that
the closed Euclidean ball of radius $r$ about $a$ is contained in $O$.
Set $V=\chi^{-1}(\{u:\|u-a\|<r\})$.
The inverse image of the corresponding closed ball is compact in $X$.
It is closed because $X$ is Hausdorff, and contains $V$, so $\overline V$
is compact.
The map
\[
 \theta:V\to\mathbb R^n,\qquad
 \theta(x)=\frac{\chi(x)-a}{r-\|\chi(x)-a\|}
\]
is a homeomorphism: the inverse sends $y\in\mathbb R^n$ to
$\chi^{-1}(a+ry/(1+\|y\|))$.
For $1\le j\le n$, define
\[
 f_j:X\to\mathbb R,\qquad
 f_j(x)=
 \begin{cases}
 \dfrac{2\theta_j(x)}{1+\|\theta(x)\|^2},&x\in V,\\[5pt]
 0,&x\notin V.
 \end{cases}
\]

We verify continuity at points outside $V$.
For $R>0$, the set
$K_R=\theta^{-1}(\{u:\|u\|\le R\})$ is a compact subset of $V$,
and hence is closed in $X$.
If $x\in V\setminus K_R$, then
\[
 |f_j(x)|\le\frac{2\|\theta(x)\|}{1+\|\theta(x)\|^2}
           \le\frac{2}{\|\theta(x)\|}<\frac2R.
\]
Given $\varepsilon>0$, choose $R>2/\varepsilon$.
The open set $X\setminus K_R$ contains $X\setminus V$, and
$|f_j|<\varepsilon$ throughout it. Since $f_j=0$ on $X\setminus V$,
this proves continuity there. Continuity in $V$ follows from the formula.
Moreover, $f_j$ vanishes outside $\overline V$.
For $\varepsilon>0$, the set $\{x:|f_j(x)|\ge\varepsilon\}$ is a
closed subset of compact $\overline V$; hence $f_j\in\C_0(X)$.

Fix $b=(b_1,\ldots,b_n)\ne0$ and put $f_b=\sum_{j=1}^n b_jf_j$.
For $x\in V$, write $u=\theta(x)$. Then
\begin{equation}\label{eq:chartmaximum}
 \|b\|-f_b(x)
 =\|b\|-\frac{2\langle b,u\rangle}{1+\|u\|^2}
 =\frac{\|b\|}{1+\|u\|^2}
   \left\|u-\frac{b}{\|b\|}\right\|^2\ge0.
\end{equation}
Equality holds precisely when $x=\theta^{-1}(b/\|b\|)$.
Outside $V$ we have $f_b(x)=0<\|b\|$, so this is the unique
point of maximum on $X$. Moreover,
\[
 |f_b(x)|\leq\|b\|\frac{2\|u\|}{1+\|u\|^2}
 \leq\|b\|\quad(x\in V),\qquad
 \|f_b\|_\infty=\|b\|.
\]
Thus $Tb=f_b$ defines the asserted linear isometry.
Every $Tb$ vanishes outside the compact set $H=\overline V$.
\end{proof}

The chart hypothesis in Lemma~\ref{lem:chart} is not necessary:
$S^{n-1}$ supports an admissible
copy of $\ell_2^n$ by the coordinate functions, although it has no
nonempty open subset homeomorphic to an open subset of $\mathbb R^n$.
For $n\ge2$ this follows from invariance of dimension, and for $n=1$
from the discreteness of $S^0$.
The guaranteed dimension $n$ cannot be increased in general:
for $n\ge2$, Theorem~\ref{thm:B} gives the exact bound $n$ on
$\mathbb R^n$; for $n=1$, Corollary~\ref{cor:compact-metric-sphere}
gives the exact bound one on $[0,1]$: it contains $S^0$, but no
subset homeomorphic to $S^1$, since every nondegenerate compact
connected subset of an interval has a cut point, whereas $S^1$ has none.

On a countably compact manifold, the evaluation map yields an embedded
sphere. We use invariance of domain in the following form: every
continuous injection from an open subset of $\mathbb R^n$ into
$\mathbb R^n$ has open image. Consequently, a continuous injective
map between topological $n$-manifolds without boundary is open when
restricted to coordinate neighborhoods
\cite[Theorem 2B.3]{Hatcher}.

\begin{lemma}\label{lem:compactdimension}
Let $M$ be a connected countably compact Hausdorff space locally
homeomorphic to $\mathbb R^n$, where $n\ge1$. No second-countability
assumption is imposed in this lemma.
\begin{enumerate}
\item Every linear subspace of $\Chat(M)\cup\{0\}$ has dimension at most
$n+1$.
\item The set $\Chat(M)$ is $(n+1)$-lineable if and only if
$M\cong S^n$. If $M\not\cong S^n$, its largest admissible dimension is $n$.
\end{enumerate}
\end{lemma}

\begin{proof}
\textup{(1)} Suppose an admissible subspace had dimension at least $n+2$.
Selecting $n+2$ linearly independent members and applying
Lemmas~\ref{lem:countable-evaluation} and~\ref{lem:gauge} would give
an embedding $e:S^{n+1}\to M$.
Choose a coordinate neighborhood $U$ in $S^{n+1}$ and a coordinate
neighborhood $V$ in $M$ such that $e(U)\subset V$.
Let $\chi:U\to O\subset\mathbb R^{n+1}$ and
$\psi:V\to P\subset\mathbb R^n$ be the corresponding charts.
Then
\[
 h=\psi\circ e\circ\chi^{-1}:O\to\mathbb R^n
\]
is a continuous injection from a nonempty open subset of
$\mathbb R^{n+1}$.
Composing it with $i:\mathbb R^n\to\mathbb R^{n+1}$,
$i(u)=(u,0)$, gives a continuous injection $i\circ h:O\to\mathbb R^{n+1}$.
Invariance of domain \cite[Theorem 2B.3]{Hatcher} makes its image open
in $\mathbb R^{n+1}$. But that image lies in
$\mathbb R^n\times\{0\}$, which has empty interior, a contradiction.
This also excludes infinite-dimensional admissible subspaces, since
any such subspace contains $n+2$ independent members.

\smallskip
\noindent\textup{(2)} If $\Chat(M)$ is $(n+1)$-lineable,
Lemmas~\ref{lem:countable-evaluation} and~\ref{lem:gauge} give
an embedding $e:S^n\to M$.
In coordinate charts, invariance of domain makes $e$ an open map.
Thus $e(S^n)$ is open in $M$. It is also closed, since it is compact
and $M$ is Hausdorff. Connectedness of $M$ and nonemptiness of
$e(S^n)$ therefore give $e(S^n)=M$.
Thus $M\cong S^n$.
Conversely, if $M\cong S^n$, the zero-set $M$ itself, defined by the
constant zero function, satisfies Proposition~\ref{prop:countably-compact-sphere} with $k=n+1$.
Hence $\Chat(M)$ is $(n+1)$-lineable.
If $M$ is not homeomorphic to $S^n$, part \textup{(1)} and the exclusion of dimension
$n+1$ give the upper bound $n$; Lemma~\ref{lem:chart} gives the matching
lower bound.
\end{proof}

\subsection{The upper bound on noncompact manifolds}

A second-countable countably compact space is compact, so the
preceding lemma leaves the noncompact case open. For functions in
$\C_0(M)$, extension by zero to the one-point compactification may
introduce a second maximizer. The boundary case of
Lemma~\ref{lem:basis} accounts for this obstruction.

\begin{lemma}\label{lem:noncompactdimension}
Let $M$ be a connected noncompact Hausdorff second-countable
$n$-manifold without boundary, where $n\ge2$.
Every linear subspace of $\Chat_0(M)\cup\{0\}$ has dimension at most $n$.
\end{lemma}

\begin{proof}
Suppose, towards a contradiction, that there exist linearly
independent functions $f_1,\ldots,f_{n+1}\in\C_0(M)$ whose span
is contained in $\Chat_0(M)\cup\{0\}$.
Let $L=M\cup\{\infty\}$ be the one-point compactification of $M$,
extend each $f_j$ continuously to $L$ by setting $f_j(\infty)=0$,
and define
\[
 F:L\longrightarrow\mathbb R^{n+1},
 \qquad F(x)=(f_1(x),\ldots,f_{n+1}(x)),
 \qquad C=\conv F(L).
\]
By Lemma~\ref{lem:basis}\textup{(1)}, the set $C$ is compact,
$0=F(\infty)\in C$, and $\interior C\ne\varnothing$.
Lemma~\ref{lem:gauge}\textup{(2)} therefore gives
$\partial C\cong S^n$.

Suppose, for a contradiction, that $0\in\interior C$.
In this case, Lemma~\ref{lem:basis}\textup{(2)} gives a
homeomorphism
$F|_{F^{-1}(\partial C)}:F^{-1}(\partial C)\to\partial C$.
Since $F(\infty)=0\notin\partial C$, its inverse defines an
embedding of $S^n$ into $M$.
By invariance of domain in coordinate charts, the image of this
embedding is open in $M$. It is also compact and hence closed,
because $M$ is Hausdorff.
Connectedness of $M$ then implies that the image is all of $M$,
contradicting noncompactness. Thus $0\in\partial C$.

Lemma~\ref{lem:basis}\textup{(3)} now provides $x_0\in M$ such
that, writing $v_0=F(x_0)$ and $A=[0,v_0]\subset\partial C$,
every point of $\partial C\setminus A$ has exactly one preimage
in $M$, and
\[
 F^{-1}(A)=\{x_0,\infty\}.
\]
We claim that $v_0=0$.
The set $A$ has empty interior relative to $\partial C$. This is immediate if $v_0=0$.
If $v_0\ne0$ and $A$ contained a nonempty relatively open set,
a coordinate chart of the $n$-manifold $\partial C$ inside that
set would give a continuous injection from a nonempty open
subset of $\mathbb R^n$ into the line $\mathbb Rv_0$.
Identifying this line with $\mathbb R$ and then including it
as a coordinate line in $\mathbb R^n$, invariance of domain
would make the resulting image open in $\mathbb R^n$.
This is impossible, since $n\ge2$ and a line has empty interior
in $\mathbb R^n$.

Since $A$ is closed and has empty relative interior,
$\partial C\setminus A$ is dense in $\partial C$.
Moreover, $F(L)$ is compact, hence closed in $\mathbb R^{n+1}$,
and contains $\partial C\setminus A$. Therefore
$\partial C\subset F(L)$.
If $v_0\ne0$, the point $v_0/2\in A$ would consequently have a
preimage in $L$. Such a preimage would belong to
$F^{-1}(A)=\{x_0,\infty\}$, although
$F(x_0)=v_0\ne v_0/2$ and $F(\infty)=0\ne v_0/2$.
This contradiction proves
\begin{equation}\label{eq:collapsedsegment}
 F(x_0)=0,\qquad A=\{0\},\qquad
 F^{-1}(\{0\})=\{x_0,\infty\}.
\end{equation}

For each $u\in\partial C\setminus\{0\}$, denote its unique
preimage by $s(u)$. This defines an injective map
\[
 s:\partial C\setminus\{0\}\longrightarrow M\setminus\{x_0\},
 \qquad F(s(u))=u.
\]
We verify that $s$ is continuous.
Every point of $\partial C\setminus\{0\}$ has a compact
neighborhood $H$ contained in this punctured boundary.
The set $F^{-1}(H)$ is compact in $L$ and avoids both $x_0$
and $\infty$.
The restriction $F:F^{-1}(H)\to H$ is a continuous bijection
from a compact space to a Hausdorff space, so its inverse
$s|_H$ is continuous.
Since the interiors of these neighborhoods cover
$\partial C\setminus\{0\}$, the map $s$ is continuous.

Put $U=s(\partial C\setminus\{0\})$.
Because $\partial C\cong S^n$, we have
$\partial C\setminus\{0\}\cong\mathbb R^n$.
Invariance of domain in coordinate charts shows that $U$ is
open in $M$ and that $s$ is a homeomorphism onto $U$.
Furthermore,
\begin{equation}\label{eq:boundaryU}
 \overline U^{\,M}\setminus U\subset\{x_0\}.
\end{equation}
Indeed, if $x\in\overline U^{\,M}\setminus U$, continuity of
$F$ and closedness of $\partial C$ imply $F(x)\in\partial C$.
If $F(x)\ne0$, uniqueness of its preimage would give
$x=s(F(x))\in U$, a contradiction.
Thus $F(x)=0$, and \eqref{eq:collapsedsegment} gives $x=x_0$.
It follows that $U$ is both open and closed in
$M\setminus\{x_0\}$.

The space $M\setminus\{x_0\}$ is connected.
Its connected components are open, since it is locally path
connected, and closed in $M\setminus\{x_0\}$.
If a component $V$ did not have $x_0$ in its closure in $M$,
then $V$ would also be closed in $M$.
As $V$ is open in $M$, this would contradict connectedness of
$M$.
Thus every component has $x_0$ in its closure.
Choose an open coordinate ball $B_0$ centered at $x_0$.
The set $B_0\setminus\{x_0\}$ is path connected because
it is homeomorphic to a punctured open ball in $\mathbb R^n$ and
$n\ge2$. It meets every component of $M\setminus\{x_0\}$.
Since a connected set lies in a single connected component,
there can be only one such component.
Consequently, the nonempty relatively open and closed set
$U$ equals $M\setminus\{x_0\}$, and hence
\begin{equation}\label{eq:puncturedEuclidean}
 M\setminus\{x_0\}=U\cong\mathbb R^n.
\end{equation}

Finally, we show directly that this forces $M$ to be compact.
Let $h:M\setminus\{x_0\}\to\mathbb R^n$ be a homeomorphism.
Choose a compact coordinate ball $B$ with $x_0$ in its interior,
and put $V=\operatorname{int}_M B$.
Since $\partial_M B=B\setminus V$ is compact and avoids $x_0$,
there exists $R>0$ such that
$h(\partial_M B)\subset\{y:\|y\|<R\}$.
Define
\[
 D=\{y\in\mathbb R^n:\|y\|>R\},
 \qquad
 K=h^{-1}\bigl(\{y\in\mathbb R^n:\|y\|\le R\}\bigr).
\]
The set $K$ is compact in $M$, hence closed, and does not
contain $x_0$.
Thus $V\setminus K$ is an open neighborhood of $x_0$.
Since $x_0$ is not isolated, this neighborhood contains a point
$x\ne x_0$, for which
$h(x)\in D\cap h(V\setminus\{x_0\})$.

The set $D$ is path connected because $n\ge2$: its radial
projection onto the sphere of radius $R+1$ deforms it to a
path-connected sphere. It is also disjoint from $h(\partial_M B)$. Hence
\[
 D=
 \bigl(D\cap h(V\setminus\{x_0\})\bigr)
 \cup
 \bigl(D\cap h(M\setminus B)\bigr)
\]
is a union of two disjoint relatively open sets.
The first is nonempty, so the second is empty.
Therefore $M\setminus B\subset K$, and consequently $M=B\cup K$.
Both $B$ and $K$ are compact, making $M$ compact.
This contradicts the hypothesis.

We have excluded admissible subspaces of dimension $n+1$.
Any subspace of larger finite or infinite dimension contains
$n+1$ linearly independent functions and is therefore excluded
as well. Thus every linear subspace of
$\Chat_0(M)\cup\{0\}$ has dimension at most $n$.
The restriction $n\ge2$ was used to show that $A$ has empty relative
interior in $\partial C$, that a coordinate ball punctured at its center
is connected, and that the exterior of a Euclidean ball is connected.
These are precisely the steps unavailable in dimension one.
\end{proof}

\begin{proof}[Proof of Theorem~\ref{thm:B}]
If $M\cong S^n$, Lemma~\ref{lem:compactdimension} gives largest admissible
dimension $n+1$. For compact $M\not\cong S^n$, the same lemma gives
dimension $n$. For noncompact $M$, Lemma~\ref{lem:chart} supplies
dimension $n$, and Lemma~\ref{lem:noncompactdimension} excludes larger
dimensions.

Thus the possible positive dimensions are precisely those stated in
the theorem, and all are finite. Theorem~\ref{thm:A} gives both extension
properties exactly when $\alpha=\beta$ lies in the corresponding range.
\end{proof}

\section{Spaceability on countably compact domains}\label{sec:countable-spaceability}

Closed infinite-dimensional admissible subspaces can exist even when
every prescribed infinite-dimensional spaceability property fails.
The obstruction comes from a sequence converging to a nonzero constant.

\subsection{Infinite prescribed dimensions}

For the supremum norm, the following construction is the obstruction
from \cite[Theorem 2.1]{FPRR}. Its uniform convergence also gives the
compact-open conclusion.

\begin{lemma}\label{lem:constantshift}
Let $X$ be a topological space with at least two points, and let
$\alpha\ge\aleph_0$ be a cardinal.
If $\Chat(X)$ is $\alpha$-lineable, there is an $\alpha$-dimensional
subspace $V\subset\Chat(X)\cup\{0\}$ for which no closed subspace of
$\C(X)$ containing $V$ is contained in $\Chat(X)\cup\{0\}$.
The same conclusion holds for closed extensions in $\C_b(X)$.
Consequently, for every $\aleph_0\leq\alpha\leq\beta$, the set
$\Chat(X)$ is not $(\alpha,\beta)$-spaceable in either ambient topology.
\end{lemma}

\begin{proof}
Choose an $\alpha$-dimensional subspace
$E\subset\Chat(X)\cup\{0\}$ and a Hamel basis $(f_i)_{i\in I}$,
where $\card(I)=\alpha$.
Every $f_i$ is nonzero and bounded, so replacing it by
$f_i/\|f_i\|_\infty$ gives a basis with $\|f_i\|_\infty=1$ for all $i$.
Let $\one:X\to\mathbb R$ be the constant function with value one.
A nonzero constant has at least two maximizing points on $X$, so
\begin{equation}\label{eq:noconstants}
 E\cap\spanop\{\one\}=\{0\}.
\end{equation}
Since $I$ is infinite, choose distinct indices $i_m\in I$ for $m\in\mathbb N$.
Define
\[
 t_i=
 \begin{cases}
 1/m,&i=i_m\text{ for some }m\in\mathbb N,\\
 1,&i\notin\{i_m:m\in\mathbb N\},
 \end{cases}
 \qquad
 \varphi_i:X\to\mathbb R,\qquad \varphi_i=\one+t_if_i,
\]
and put $V=\spanop\{\varphi_i:i\in I\}$.

We have $\dim V=\alpha$.
Let $J\subset I$ be finite, and suppose $\sum_{i\in J}a_i\varphi_i=0$.
Expanding the definition gives
\[
 \sum_{i\in J}a_it_if_i
 =-\left(\sum_{i\in J}a_i\right)\one.
\]
The left side lies in $E$ and the right side in $\spanop\{\one\}$.
By \eqref{eq:noconstants}, both sides are zero.
Independence of the $f_i$ gives $a_it_i=0$ for every $i\in J$.
Since $t_i>0$, all coefficients $a_i$ vanish.
Thus $(\varphi_i)_{i\in I}$ is independent, and $\dim V=\card(I)=\alpha$.

Next, take a nontrivial finite coefficient family $(a_i)_{i\in J}$ and put
\[
 f=\sum_{i\in J}a_it_if_i,\qquad c=\sum_{i\in J}a_i.
\]
Independence of the $f_i$ and positivity of $t_i$ give $f\ne0$.
The corresponding element of $V$ is $f+c\one$.
For every $x\in X$,
\[
 \begin{aligned}
 (f+c\one)(x)=\sup_{y\in X}(f(y)+c)
 &\quad\Longleftrightarrow\quad f(x)+c=\sup_{y\in X}f(y)+c\\
 &\quad\Longleftrightarrow\quad f(x)=\sup_{y\in X}f(y).
 \end{aligned}
\]
Hence $\Argmax_X(f+c\one)=\Argmax_Xf$, which is a singleton.
Thus $V\subset\Chat(X)\cup\{0\}$.

Finally,
\begin{equation}\label{eq:constantlimit}
 \|\varphi_{i_m}-\one\|_\infty
 =t_{i_m}\|f_{i_m}\|_\infty=\frac1m\longrightarrow0.
\end{equation}
For each nonempty compact $H\subset X$, the defining seminorm of the
compact-open topology satisfies
\[
 \sup_{x\in H}|\varphi_{i_m}(x)-1|
 \le\|\varphi_{i_m}-\one\|_\infty=\frac1m.
\]
Thus $\varphi_{i_m}\to\one$ both uniformly and in the compact-open topology.
Any subspace closed in either ambient topology and containing $V$ must
therefore contain $\one$.
But $\one\ne0$ and $\Argmax_X\one=X$ has at least two points.
Such a closed subspace cannot be contained in $\Chat(X)\cup\{0\}$.
If $\alpha$-lineability fails, $(\alpha,\beta)$-spaceability fails
by definition. Otherwise, the subspace just constructed excludes it
for every $\beta\geq\alpha$, proving the final assertion.
\end{proof}

The restriction to infinite $\alpha$ in this lemma is essential.
Every finite-dimensional admissible subspace is already closed in
both ambient topologies and is therefore a closed admissible extension
of itself. The compact-open topology is Hausdorff on $\C(X)$ for
any topological space $X$, since its seminorms include evaluation at
each singleton. Thus no separation axiom or local compactness is
needed for the infinite-dimensional obstruction.

\begin{corollary}\label{cor:countable-classification}
Let $X$ be a countably compact space with at least two points whose
continuous real functions separate points, and let $1\leq\alpha\leq\beta$
be cardinals. Then $\Chat(X)$ is $(\alpha,\beta)$-spaceable if and only
if there is $k\in\mathbb N$ such that $\alpha=\beta=k$ and $X$ contains
a zero-set homeomorphic to $S^{k-1}$.
The equivalence holds in $\C_b(X)$ with the supremum norm and in
$\C(X)$ with the compact-open topology.
\end{corollary}

\begin{proof}
The hypothesis on continuous functions implies that $X$ is Hausdorff.
For finite $\alpha$, the countably compact case of Theorem~\ref{thm:A}
says that $\Chat(X)$ is $(\alpha,\beta)$-spaceable precisely when
$\alpha=\beta$ and $\Chat(X)$ is $\alpha$-lineable, in either topology.
Writing $k=\alpha$, Proposition~\ref{prop:countably-compact-sphere}
identifies this last condition with the existence of the required zero-set.
For infinite $\alpha$, Lemma~\ref{lem:constantshift} excludes
$(\alpha,\beta)$-spaceability in both topologies.
\end{proof}

\subsection{Closed admissible subspaces from Banach spaces}\label{sec:example}

The Hilbert-space example of \cite[Example 4.1]{BCFPS} extends to
reflexive strictly convex Banach spaces: the weak unit ball is compact,
and every nonzero continuous linear functional has a unique maximizer.
Recall that a norm is strictly convex if
$\|(x+y)/2\|<1$ whenever $x\ne y$ and $\|x\|=\|y\|=1$.

\begin{proposition}\label{prop:dual-realization}
Let $E$ be a nonzero real Banach space with continuous dual $E^*$,
and put $B_E=\{x\in E:\|x\|\le1\}$. Let
$K=(B_E,\sigma(E,E^*))$, where $\sigma(E,E^*)$ is the weak topology,
that is, the coarsest topology on $E$ making every functional in
$E^*$ continuous, restricted to $B_E$. The following are equivalent:
\begin{enumerate}
\item $E$ is reflexive and strictly convex.
\item $K$ is compact and every nonzero $f\in E^*$ has a unique
point of maximum on $K$.
\end{enumerate}
When these conditions hold, $K$ is compact Hausdorff and
\[
 J:E^*\longrightarrow\C(K),\qquad (Jf)(x)=f(x),
\]
is a linear isometry with closed image satisfying
$J(E^*)\setminus\{0\}\subset\Chat(K)$.
If $E$ is separable, then $K$ is metrizable.
\end{proposition}

\begin{proof}
The weak topology is Hausdorff because continuous linear functionals
separate points by the Hahn--Banach theorem. Kakutani's theorem
\cite[Theorem 3.17]{Brezis} states that $B_E$ is weakly compact if
and only if $E$ is reflexive.

Assume (1). Every $f\in E^*$ is continuous on $K$ by the definition
of the weak topology. Compactness and symmetry of $B_E$ imply
\[
 \max_{x\in B_E}f(x)=\sup_{x\in B_E}|f(x)|=\|f\|.
\]
If $f\ne0$ and $f(x)=\|f\|$ for some $x\in B_E$, then
$\|f\|=f(x)\leq\|f\|\|x\|$ forces $\|x\|=1$.
If distinct $x,y\in B_E$ both maximized $f$, strict convexity
would give
\[
 \|f\|=f\left(\frac{x+y}{2}\right)
 \leq\|f\|\left\|\frac{x+y}{2}\right\|<\|f\|,
\]
a contradiction. This proves (2).

Conversely, (2) and Kakutani's theorem give reflexivity.
If the norm were not strictly convex, there would exist distinct
$x,y\in E$ with $\|x\|=\|y\|=\|(x+y)/2\|=1$.
By Hahn--Banach, some $f\in E^*$ satisfies
$\|f\|=1$ and $f((x+y)/2)=1$. Since $f(x),f(y)\leq1$ and
$(f(x)+f(y))/2=1$, both $f(x)$ and $f(y)$ equal one.
They are distinct maximizing points, contradicting (2).
Thus (1) follows.

The formula defining $J$ is linear, and
$\|Jf\|_\infty=\sup_{x\in B_E}|f(x)|=\|f\|$.
The image is closed because $E^*$ is complete and $J$ is an isometry.
Finally, if $E$ is separable and reflexive, then $E^*$ is separable
\cite[Corollary 3.27]{Brezis}, and the weak topology on $B_E$ is
metrizable by \cite[Theorem 3.29]{Brezis}.
\end{proof}

Both assumptions are necessary for the weak unit-ball
construction. For example, $E=(\mathbb R^2,\|\cdot\|_\infty)$ is
reflexive but not strictly convex, and $f(x_1,x_2)=x_1$ is maximized
at every point of $\{1\}\times[-1,1]$.

For $1\le r<\infty$, let $\ell_r$ denote the Banach space of real
sequences $a=(a_j)_{j\in\N}$ satisfying
$\sum_{j=1}^{\infty}|a_j|^r<\infty$, with norm
$\|a\|_r=(\sum_{j=1}^{\infty}|a_j|^r)^{1/r}$.
For $E=\ell_q$, this gives compact metric domains with closed
admissible copies of $\ell_p$.

\begin{corollary}\label{cor:lp-realization}
For each $1<p<\infty$, put $q=p/(p-1)$ and
$K_q=(B_{\ell_q},\sigma(\ell_q,\ell_p))$.
Then $K_q$ is compact metrizable, and
\[
 J_p:\ell_p\longrightarrow\C(K_q),\qquad
 (J_pa)(x)=\sum_{j=1}^{\infty}a_jx_j,
\]
is a linear isometry with closed admissible image.
In particular, $\Chat(K_q)$ is maximally spaceable.
For $1\leq\alpha<\aleph_0$ and $\beta\geq\alpha$,
\[
 \Chat(K_q)\text{ is }(\alpha,\beta)\text{-lineable}
 \quad\Longleftrightarrow\quad\alpha=\beta.
\]
For all cardinals $1\leq\alpha\leq\beta$,
\[
 \Chat(K_q)\text{ is }(\alpha,\beta)\text{-spaceable}
 \quad\Longleftrightarrow\quad\alpha=\beta<\aleph_0.
\]
For $a\ne0$, the unique maximizing point is
\[
 x_j^a=\frac{\operatorname{sgn}(a_j)|a_j|^{p-1}}{\|a\|_p^{p-1}}
 \quad(j\in\mathbb N),
\]
where $\operatorname{sgn}(0)=0$.
\end{corollary}

\begin{proof}
The space $\ell_q$ is separable, reflexive and strictly convex,
and its dual is isometrically $\ell_p$. These standard properties
of the $L^p$ spaces apply to counting measure; see \cite[Chapter 4]{Brezis}.
Proposition~\ref{prop:dual-realization} therefore gives compactness,
metrizability and the asserted isometry. H\"older's inequality gives
absolute convergence of the defining series. Since $(p-1)q=p$,
\[
 \sum_{j=1}^{\infty}|x_j^a|^q
 =\frac{\sum_{j=1}^{\infty}|a_j|^p}{\|a\|_p^p}=1,
 \qquad
 (J_pa)(x^a)=\frac{\sum_{j=1}^{\infty}|a_j|^p}{\|a\|_p^{p-1}}
 =\|a\|_p.
\]
The proposition ensures uniqueness, so this is the maximizing point.

We verify the dimension assertion explicitly.
The cardinality of $\ell_p\subset\mathbb R^{\mathbb N}$ is at most
$\cc$. For every $0<t<1$, the sequence $a^t=(t^j)_{j\geq1}$ lies
in $\ell_p$, since $\sum_{j\geq1}t^{jp}=t^p/(1-t^p)<\infty$.
For distinct $t_1,\ldots,t_r\in(0,1)$, a relation
$\sum_{i=1}^r c_i a^{t_i}=0$ gives
$\sum_{i=1}^r c_it_i^j=0$ for $1\leq j\leq r$.
The determinant
\[
 \det(t_i^j)_{j,i=1}^r
 =\left(\prod_{i=1}^r t_i\right)
   \prod_{1\leq i<\ell\leq r}(t_\ell-t_i)
\]
is nonzero. Hence all $c_i$ vanish and $\dim\ell_p=\cc$.
Choose a countable dense subset $D$ of the compact metric space $K_q$.
Continuous functions are determined by their restrictions to $D$, so
\[
 \cc=\dim J_p(\ell_p)\leq\dim\C(K_q)
 \leq\card\C(K_q)\leq\card(\mathbb R^D)\leq\cc.
\]
Thus the closed image is maximally dimensional.
It also supplies admissible subspaces of every positive finite
dimension. Since $K_q$ is compact, Theorem~\ref{thm:A} yields the
finite-dimensional extension classification.
Lemma~\ref{lem:constantshift} excludes every infinite prescribed
dimension for spaceability, proving the remaining assertion.
\end{proof}

The range $1<p<\infty$ is sharp for the specific dual-ball model
in this corollary. At $p=1$, the corresponding domain would be
$B_{\ell_\infty}$ with its weak topology; it is not compact, and
the coordinate functional $x\mapsto x_1$ has many maximizers.
Even changing to the weak-star topology, which makes that ball
compact, does not restore uniqueness. At $p=\infty$, the domain
would be $B_{\ell_1}$, which is not weakly compact. Moreover,
$a=(1-1/(j+1))_{j\geq1}\in\ell_\infty$ defines a functional with
supremum one on $B_{\ell_1}$ and no maximizing point: for $x\ne0$,
\[
 \sum_{j=1}^{\infty}a_jx_j
 \leq\sum_{j=1}^{\infty}a_j|x_j|
 <\sum_{j=1}^{\infty}|x_j|\leq1,
 \qquad \langle a,e_j\rangle\longrightarrow1,
\]
while the value at $x=0$ is zero.

\section{Proof of Theorem C}\label{sec:proofC}

Each direction on the coefficient sphere determines the unique
maximizer of the corresponding function. Directions with the same
maximizer form a closed equivalence relation by boundedness of the
evaluation vectors and uniqueness of the maximizing points. Its equivalence classes are contractible.
The Vietoris--Begle theorem identifies the cohomology of this
quotient with that of the sphere; the covering dimension of the
domain bounds the dimension of the quotient.

\subsection{Topological definitions and the mapping theorem}

A space is \emph{$\sigma$-compact} if it is a countable union of
compact subsets. By replacing the first $r$ compact sets by their
union, we may always arrange a representation $X=\bigcup_{r\ge1}K_r$
with $K_r\subset K_{r+1}$ compact.
For a metrizable space $Y$, write $\dim_{\mathrm{top}}Y$ for its
\emph{covering dimension}. The inequality $\dim_{\mathrm{top}}Y\le m$
means that every finite open cover has an open refinement of
multiplicity at most $m+1$: no point belongs to more than $m+1$
members of the refinement. Here a cover $\mathcal V$ refines
$\mathcal U$ if each member of $\mathcal V$ is contained in a member
of $\mathcal U$. This is the usual covering dimension
\cite[Section 3.1]{Engelking}; the subscript distinguishes it from
Hamel dimension.
We use monotonicity on metrizable subspaces and the countable closed
sum theorem \cite[Theorems 4.1.3, 4.1.7 and 4.1.9]{Engelking}:
\begin{equation}\label{eq:dimensionsum}
 Y=\bigcup_{r\ge1}Y_r,\quad Y_r\text{ closed in }Y,\quad
 \dim_{\mathrm{top}}Y_r\le m
 \quad\Longrightarrow\quad\dim_{\mathrm{top}}Y\le m.
\end{equation}
Also $\dim_{\mathrm{top}}\mathbb R^n=n$
\cite[Theorem 1.8.2]{Engelking}.

We use \v{C}ech cohomology with integer coefficients. For compact
metrizable $Y$, it can be computed as
\begin{equation}\label{eq:cechdefinition}
 \check H^j(Y;\mathbb Z)
 =\varinjlim_{\mathcal U}H^j(N(\mathcal U);\mathbb Z),
\end{equation}
where $\mathcal U$ ranges over finite open covers ordered by refinement
and $N(\mathcal U)$ is the nerve of the cover. Its vertices are the
members of $\mathcal U$, and a finite family spans a simplex exactly
when its intersection is nonempty. Refinements induce the homomorphisms
in this direct limit; continuous maps induce cohomology homomorphisms
contravariantly. See \cite[Sections 6.7--6.9]{Spanier}.
A nonempty compact contractible space is acyclic, meaning that its
reduced \v{C}ech cohomology groups vanish in every degree.

We use the following precise form of the Vietoris--Begle mapping
theorem \cite[Theorem 6.9.15, p.~344]{Spanier}: if $P,Q$ are compact
metrizable spaces and $q:P\to Q$ is continuous and onto, with every
fiber nonempty and contractible, then
\begin{equation}\label{eq:vietorisbegle}
 q^*:\check H^j(Q;\mathbb Z)\longrightarrow
 \check H^j(P;\mathbb Z)
 \quad\text{is an isomorphism for every }j\ge0.
\end{equation}
Indeed, $q$ is closed by compactness, and contractibility makes its
fibers acyclic in every degree, so the finite-degree theorem applies
with an arbitrarily large degree cutoff. Spanier formulates the theorem in Alexander cohomology,
which agrees with \v{C}ech cohomology for compact metrizable spaces.

Two computations will be used. For $r\ge1$,
\begin{equation}\label{eq:cechsphere}
 \check H^r(S^r;\mathbb Z)\cong\mathbb Z.
\end{equation}
This follows from the agreement with simplicial cohomology on finite
polyhedra and the usual sphere computation; see \cite[Section 3.1]{Hatcher}.
On the other hand, for compact metrizable $Y$,
\begin{equation}\label{eq:cechdimensionvanish}
 \dim_{\mathrm{top}}Y\le m
 \quad\Longrightarrow\quad
 \check H^j(Y;\mathbb Z)=0\quad(j>m).
\end{equation}
To verify this directly, refine any finite open cover to one of
multiplicity at most $m+1$ and take a finite subcover of the
refinement. Its nerve has no simplices of dimension greater than
$m$, so its cochain groups, and hence its cohomology groups, vanish
in degrees greater than $m$. Every class from the original cover
therefore becomes zero in such a refinement, proving
\eqref{eq:cechdimensionvanish} by the definition of the direct limit.

\subsection{An upper bound on \texorpdfstring{$\sigma$}{sigma}-compact domains}

\begin{proposition}\label{prop:sigma-upper}
Let $X$ be a $\sigma$-compact metrizable space with at least two points
and $\dim_{\mathrm{top}}X\le m<\infty$, where $m$ is a nonnegative
integer. Every admissible linear subspace of $\C(X)$ has Hamel
dimension at most $m+1$.
\end{proposition}

\begin{proof}
It suffices to bound a finite-dimensional admissible subspace
$E=\operatorname{span}\{f_1,\ldots,f_d\}$ with $d\ge2$;
if a larger finite or infinite dimension existed, it would contain
such a subspace with $d=m+2$.
For $a\in\mathbb R^d$, put
\[
 f_a=\sum_{i=1}^d a_if_i,\qquad
 F(x)=(f_1(x),\ldots,f_d(x)),\qquad
 H(a)=\sup_{x\in X}\langle a,F(x)\rangle.
\]
Every $f_i$ is bounded, since both $f_i$ and $-f_i$ attain a finite
maximum. Consequently $F(X)$ is bounded and $H$ is finite, with
\begin{equation}\label{eq:Hspherecontinuity}
 |H(a)-H(b)|\le\sup_{x\in X}|f_a(x)-f_b(x)|
 \le\sum_{i=1}^d|a_i-b_i|\,\|f_i\|_\infty.
\end{equation}
In particular, $H$ is continuous. It also satisfies
$H(ta)=tH(a)$ for $t\ge0$ and
$H(a+b)\le H(a)+H(b)$, by taking suprema in the corresponding
pointwise identities and inequality.
For $a\ne0$, let $x_a$ be the unique maximizer of $f_a$.
Linear independence ensures $f_a\ne0$.
No nonzero constant belongs to $E$, since $X$ has at least two points.
We will construct a compact metrizable quotient $Q$ of $S^{d-1}$ such that
\[
 \check H^{d-1}(Q;\mathbb Z)\cong\mathbb Z,
 \qquad \dim_{\mathrm{top}}Q\le m.
\]
By \eqref{eq:cechdimensionvanish}, these properties force $d-1\le m$.

On $S=S^{d-1}$ define $a\sim b$ if $x_a=x_b$, and let
$R=\{(a,b)\in S\times S:a\sim b\}$.
The relation $R$ is closed.
Suppose $(a_j,b_j)\in R$ and $(a_j,b_j)\to(a,b)$.
Set $x_j=x_{a_j}=x_{b_j}$. Boundedness of $F(X)$ permits a
subsequence for which $F(x_j)\to v\in\mathbb R^d$.
Although $v$ need not belong to $F(X)$, continuity of $H$ gives
\begin{equation}\label{eq:limitcommonmax}
 \langle a,v\rangle=\lim_jH(a_j)=H(a),\qquad
 \langle b,v\rangle=\lim_jH(b_j)=H(b).
\end{equation}
We have $a+b\ne0$. Otherwise $b=-a$ and
$H(a)+H(-a)=0$, so $\sup_X f_a=\inf_X f_a$.
Then $f_a$ would be a nonzero constant, a contradiction.
Furthermore, since
$\langle a+b,F(x_j)\rangle\le H(a+b)$, we obtain
\[
 H(a+b)\le H(a)+H(b)
 =\langle a+b,v\rangle
 =\lim_j\langle a+b,F(x_j)\rangle\le H(a+b).
\]
Thus $H(a+b)=H(a)+H(b)$.
At the maximizer $z=x_{a+b}$,
\[
 0=H(a)+H(b)-f_{a+b}(z)
   =(H(a)-f_a(z))+(H(b)-f_b(z)).
\]
Both summands are nonnegative, so both are zero.
Hence $z=x_a=x_b$, proving $(a,b)\in R$.
The sequential argument proves closedness because $S\times S$ is
metrizable.

Let $Q=S/\!\sim$ be the set of equivalence classes and let
$q:S\to Q$ send a direction to its class.
We give $Q$ the \emph{quotient topology}: $O\subset Q$ is open
exactly when $q^{-1}(O)$ is open in $S$.
We verify that $Q$ is compact metrizable.
It is compact as a continuous image of $S$.
If $A\subset S$ is closed, then
\[
 q^{-1}(q(A))
 =\{a\in S:\text{there is }b\in A\text{ with }(a,b)\in R\}
\]
is the projection of the compact set $R\cap(S\times A)$, and is
therefore closed. By the definition of the quotient topology, $q(A)$
is closed in $Q$. Thus $q$ is closed.
Each class is a closed subset of $S$, hence compact.
For two distinct classes, choose disjoint open neighborhoods $U,V$
of these compact subsets in $S$. The open sets
$Q\setminus q(S\setminus U)$ and $Q\setminus q(S\setminus V)$
contain the respective classes and are disjoint. Hence $Q$ is Hausdorff.

Take a countable base of $S$ and let $\mathcal B$ consist of its
finite unions. The countable collection
\[
 \{Q\setminus q(S\setminus U):U\in\mathcal B\}
\]
is a base of $Q$: if $\xi\in O$ with $O$ open, the compact class
$q^{-1}(\{\xi\})$ is covered by finitely many basic open subsets
contained in $q^{-1}(O)$; their union $U$ gives
$\xi\in Q\setminus q(S\setminus U)\subset O$.
Thus $Q$ is second countable. A compact Hausdorff space is regular,
so the Urysohn metrization theorem
\cite[Theorem 34.1, p.~215]{Munkres} makes $Q$ metrizable.

Each fiber of $q$ is contractible.
Fix a class $A$ and $a_0\in A$, and write $x=x_{a_0}$.
The class cannot contain antipodal points: if $a,-a\in A$, then
$x$ maximizes both $f_a$ and $-f_a$, making $f_a$ constant.
Therefore
\begin{equation}\label{eq:classcontraction}
 \Phi:A\times[0,1]\to S,\qquad
 \Phi(a,t)=\frac{(1-t)a+ta_0}{\|(1-t)a+ta_0\|}
\end{equation}
is well defined and continuous. Indeed, a zero denominator with
$0<t<1$ would give $a=-a_0$ (and $t=1/2$), which was excluded.
For $y\ne x$, both $f_a(x)-f_a(y)$ and
$f_{a_0}(x)-f_{a_0}(y)$ are strictly positive. Consequently
\[
 f_{(1-t)a+ta_0}(x)-f_{(1-t)a+ta_0}(y)
 =(1-t)(f_a(x)-f_a(y))+t(f_{a_0}(x)-f_{a_0}(y))>0.
\]
Dividing by the positive denominator in \eqref{eq:classcontraction}
shows that $x_{\Phi(a,t)}=x$. Thus $\Phi$ takes values in $A$,
with $\Phi(a,0)=a$ and $\Phi(a,1)=a_0$.
Every fiber is therefore contractible. Equations
\eqref{eq:vietorisbegle} and \eqref{eq:cechsphere} yield
\begin{equation}\label{eq:quotientcohomology}
 \check H^{d-1}(Q;\mathbb Z)
 \cong\check H^{d-1}(S^{d-1};\mathbb Z)
 \cong\mathbb Z.
\end{equation}

It remains to bound $\dim_{\mathrm{top}}Q$ from the domain $X$.
Choose an increasing compact cover $X=\bigcup_{r\ge1}K_r$ and put
\[
 T_r=\{a\in S:x_a\in K_r\},\qquad Q_r=q(T_r).
\]
The set $T_r$ is closed in $S$. To check this, let $a_j\in T_r$
converge to $a$. A subsequence satisfies $x_{a_j}\to x\in K_r$,
and continuity of $F$ and \eqref{eq:Hspherecontinuity} give
$f_a(x)=\lim_j f_{a_j}(x_{a_j})=\lim_j H(a_j)=H(a)$.
Uniqueness implies $x=x_a\in K_r$.
Hence $T_r$ is compact and $Q_r$ is compact, thus closed in $Q$.

The map $m_r:T_r\to K_r$, $m_r(a)=x_a$, is continuous.
Indeed, for any sequence $a_j\to a$ in $T_r$, every convergent
subsequence of the points $x_{a_j}$ has limit $x_a$, by the same
calculation. If $x_{a_j}$ did not converge to $x_a$, compactness of
$K_r$ would give a subsequence staying outside some neighborhood of
$x_a$ and a further convergent subsequence, a contradiction.

The map $q|_{T_r}:T_r\to Q_r$ is a continuous closed surjection,
since its domain is compact and its codomain is Hausdorff.
It is therefore a quotient map. Since $m_r$ is constant precisely
on its equivalence classes, it induces a continuous bijection
\[
 \overline m_r:Q_r\to m_r(T_r)\subset K_r,\qquad
 \overline m_r(q(a))=x_a.
\]
For injectivity, $x_a=x_b$ implies $a\sim b$ and $q(a)=q(b)$;
surjectivity follows from the definition of $m_r(T_r)$.
Its domain is compact and its codomain Hausdorff, so it is a
homeomorphism. By monotonicity,
$\dim_{\mathrm{top}}Q_r=\dim_{\mathrm{top}}m_r(T_r)\le m$.
Every maximizer belongs to some $K_r$, so $S=\bigcup_rT_r$ and
$Q=\bigcup_rQ_r$. The countable closed sum theorem
\eqref{eq:dimensionsum} now gives $\dim_{\mathrm{top}}Q\le m$.
If $d-1>m$, \eqref{eq:cechdimensionvanish} contradicts
\eqref{eq:quotientcohomology}. Therefore $d\le m+1$.
\end{proof}

The maximizing-point map need not be continuous on the whole
coefficient sphere. Restricting to the compact sets $T_r$ is what
allows the countable closed sum theorem to be applied to $Q$.

\subsection{An isometric construction attaining the Euclidean bound}

The lower-bound construction follows the sphere-extension method of
\cite[Theorem 5]{AB} and \cite[Theorem 2.1]{BCFPS}. A closed subset of $\mathbb R^n$ admitting a
continuous bijection onto $S^n$ gives the points at which the maxima
are attained.

\begin{lemma}\label{lem:sharpconstruction}
For every $n\ge1$, there is a linear isometry
$T:\ell_2^{n+1}\to\C_b(\mathbb R^n)$ whose nonzero image functions
attain their maximum at exactly one point.
\end{lemma}

\begin{proof}
We first construct, by induction, a closed proper subset
$A_n\subset\mathbb R^n$ and a continuous bijection
$u_n:A_n\to S^n$.
For $n=1$, set
\[
 A_1=[0,\infty),\qquad
 u_1(t)=\left(\cos\frac{2\pi t}{1+t},\sin\frac{2\pi t}{1+t}\right).
\]
The function $t\mapsto t/(1+t)$ maps $[0,\infty)$ bijectively onto
$[0,1)$, with inverse $s\mapsto s/(1-s)$, and the usual angular
parametrization maps $[0,1)$ bijectively onto $S^1$.
Thus $u_1$ is continuous and bijective. Only a continuous bijection is
needed: the inverse is not required to be continuous. In fact, no
$u_n$ constructed below is a homeomorphism, since $A_n$ is noncompact
whereas $S^n$ is compact.

For $n\ge2$, assume $A_{n-1}$ and $u_{n-1}$ have been constructed.
The product $P=A_{n-1}\times\mathbb R$ is closed and proper in
$\mathbb R^n$. Its complement is nonempty and open, so choose
 distinct points $p_+,p_-\notin P$ and put
\[
 A_n=P\cup\{p_+,p_-\},\qquad
 u_n(x,t)=\left(\frac{u_{n-1}(x)}{\sqrt{1+t^2}},
                   \frac{t}{\sqrt{1+t^2}}\right),\qquad
 u_n(p_\pm)=(0,\ldots,0,\pm1).
\]
The set $A_n$ is closed and remains proper, because removing two
points from the nonempty open set $\mathbb R^n\setminus P$ does not
exhaust it. The two added points are isolated in $A_n$, and $P$
is open and closed in $A_n$, so the displayed formulas define a
continuous map. On $P$,
\[
 \|u_n(x,t)\|^2
 =\frac{\|u_{n-1}(x)\|^2+t^2}{1+t^2}=1.
\]
For a point $(v,s)\in S^n$ with $|s|<1$, its unique preimage in $P$
is determined by
\[
 t=\frac{s}{\sqrt{1-s^2}},\qquad
 x=u_{n-1}^{-1}\left(\frac{v}{\sqrt{1-s^2}}\right).
\]
Indeed, $\|v\|^2=1-s^2$, and substitution gives
$t/\sqrt{1+t^2}=s$ and
$u_{n-1}(x)/\sqrt{1+t^2}=v$.
Neither pole is obtained on $P$, since $|t|/\sqrt{1+t^2}<1$.
The two added points supply precisely these poles, proving bijectivity.

Write $u_n=(u_{n,1},\ldots,u_{n,n+1})$.
The Tietze extension theorem
\cite[Theorem 35.1, p.~219]{Munkres} applies to the closed subset $A_n$ of
the normal space $\mathbb R^n$. It gives continuous
$B_j:\mathbb R^n\to[-1,1]$ extending $u_{n,j}$.
Put $B=(B_1,\ldots,B_{n+1}):\mathbb R^n\to\mathbb R^{n+1}$
and define $G:\mathbb R^n\to\mathbb R^{n+1}$ by
\begin{equation}\label{eq:sharpG}
 G(x)=\frac{B(x)}{(1+\dist(x,A_n))\max\{1,\|B(x)\|\}}
 \qquad(x\in\mathbb R^n).
\end{equation}
The distance function is continuous, since
$|\dist(x,A_n)-\dist(y,A_n)|\le\|x-y\|$.
The denominator in \eqref{eq:sharpG} is positive, so $G$ is continuous.
For $x\in A_n$, we have $B(x)=u_n(x)$, $\|B(x)\|=1$ and
$\dist(x,A_n)=0$; hence $G(x)=u_n(x)$.
For $x\notin A_n$, closedness gives $\dist(x,A_n)>0$ and
\[
 \|G(x)\|
 =\frac{\|B(x)\|}{(1+\dist(x,A_n))\max\{1,\|B(x)\|\}}
 \le\frac1{1+\dist(x,A_n)}<1.
\]
Thus $G$ takes every value of $S^n$ at exactly one point, and every
other value lies in the open unit ball.

Define $(Ta)(x)=\langle a,G(x)\rangle$ for
$a\in\mathbb R^{n+1}$. This is a linear map into $\C_b(\mathbb R^n)$,
with $|(Ta)(x)|\le\|a\|\|G(x)\|\le\|a\|$.
For $a\ne0$, put $x_a=u_n^{-1}(a/\|a\|)$.
Then $(Ta)(x_a)=\|a\|$. If $x\notin A_n$, the strict norm estimate
gives $(Ta)(x)<\|a\|$.
If $x\in A_n$, both $u_n(x)$ and $a/\|a\|$ are unit vectors, and
\[
 \left\|u_n(x)-\frac a{\|a\|}\right\|^2
 =2-\frac{2\langle a,u_n(x)\rangle}{\|a\|}.
\]
Consequently $(Ta)(x)=\|a\|$ if and only if
$u_n(x)=a/\|a\|$, equivalently $x=x_a$.
We have proved both uniqueness of the maximum and
$\|Ta\|_\infty=\|a\|$.
Thus $T$ is an isometry and its image has dimension $n+1$.
\end{proof}

\begin{proof}[Proof of Theorem~\ref{thm:C}]
Proposition~\ref{prop:sigma-upper}, applied to $\mathbb R^n$, bounds
every admissible dimension by $n+1$, while
Lemma~\ref{lem:sharpconstruction} realizes dimension $n+1$.
This proves the first assertion, including exclusion of infinite
dimensions.

Let $1\le\alpha\le\beta$ be cardinals.
If $\alpha>n+1$, ordinary $\alpha$-lineability fails, so all three
extension properties in the statement fail by definition.
If $\alpha=\beta\le n+1$, subspaces of the construction supply
$\alpha$-lineability. Every prescribed $\alpha$-dimensional admissible
subspace is its own extension and is closed in both ambient
Hausdorff topological vector spaces, proving all three properties.

Suppose $\alpha<\beta$ and $\alpha\le n+1$.
If $\alpha=n+1$, no admissible subspace of dimension $\beta$ exists,
by the upper bound, so no such extension can exist.
If $1\le\alpha\le n$, Lemma~\ref{lem:chart} supplies an
$\alpha$-dimensional admissible subspace of $\C_0(\mathbb R^n)$.
Corollary~\ref{cor:transfer} then excludes $(\alpha,\beta)$-lineability
in $\Chat(\mathbb R^n)$ and hence excludes both spaceability
properties. This argument also covers $n=1$ and $\alpha=1$.
\end{proof}

\section*{Funding}
The work of the second author was supported by the Conselho Nacional
de Desenvolvimento Cient\'ifico e Tecnol\'ogico (CNPq, Brazil),
under grants 406457/2023-9 (CNPq/MCTI No.~10/2023),
403964/2024-5 (MCTI/CNPq No.~16/2024) and 305807/2025-0.

\section*{Declarations}
\noindent\textit{Competing interests.} The authors declare no competing interests.

\medskip
\noindent\textit{Data availability.} No datasets were generated or analyzed
during the current study.

\medskip
\noindent\textit{Use of generative AI.} ChatGPT (OpenAI) was used
for exploratory calculations, consistency checks, assistance with the
organization of arguments and language editing. The mathematical
arguments and final text were reviewed by the authors, who take
responsibility for the content.

\end{document}